\input ./style/arxiv-general.cfg
\documentclass[aop,MSNbibl,secthm,seceqn,dvips]{arximspdf}
\makeatletter
   \@ifpackageloaded{graphicx}{}{\usepackage{graphicx}}
\makeatother
\usepackage{mathrsfs}

\doi{10.1214/15-AOP1045}
\volume{44}
\issue{5}
\pubyear{2016}
\firstpage{3207}
\lastpage{3233}
\docsubty{FLA}

\makeatletter
\let\leqslant\leq
\let\geqslant\geq
\newcommand{\rrvert}{\vert}
\newcommand{\rrVert}{\Vert}
\newcommand{\llvert}{\vert}
\newcommand{\llVert}{\Vert}
\newtheorem{prop}[thm]{Proposition}

\newtheorem{cor}[thm]{Corollary}
\newproclaim{rem}[thm]{Remark}
\newcommand{\tarpas}{}
\makeatother

\begin{document}
\begin{frontmatter}

\title{Correlation structure of the corrector in stochastic homogenization}
\runtitle{Correlation structure in stochastic homogenization}

\begin{aug}
\author[A]{\fnms{Jean-Christophe}~\snm{Mourrat}\corref{}\ead[label=e1]{jean-christophe.mourrat@ens-lyon.fr}}
\and
\author[B]{\fnms{Felix}~\snm{Otto}\ead[label=e2]{otto@mis.mpg.de}}
\runauthor{J.-C. Mourrat and F. Otto}
\affiliation{CNRS and Max Planck Institute for Mathematics
in the Sciences}
\address[A]{CNRS\\
Ecole normale sup\'erieure de Lyon\\
46, all\'ee d'Italie\\
69007 Lyon\\
France\\
\printead{e1}}
\address[B]{Max Planck Institute for Mathematics\\
\quad in the Sciences\\
Inselstr. 22\\
04103 Leipzig\\
Germany\\
\printead{e2}}
\end{aug}

%
\received{\smonth{2} \syear{2014}}
%
\revised{\smonth{1} \syear{2015}}

\begin{abstract}
Recently, the quantification of errors in the stochastic homogenization
of divergence-form operators has witnessed important progress. Our aim
now is to go beyond error bounds, and give precise descriptions of the
effect of the randomness, in the large-scale limit. This paper is a
first step in this direction. Our main result is to identify the
correlation structure of the corrector, in dimension $3$ and higher.
This correlation structure is similar to, but different from that of a
Gaussian free field.
\end{abstract}

\end{frontmatter}

\section{Introduction}
\label{sintro}

Consider the solution $u_\varepsilon: \mathbb{R}^d \to\mathbb{R}$
of the equation
\[
\biggl(1-\nabla\cdot A \biggl(\frac{\cdot}{\varepsilon} \biggr) \nabla \biggr)
u_\varepsilon= f,
\]
where $f$ is a bounded smooth function, $A$ is a random field of
symmetric matrices on $\mathbb{R}^d$, and $\varepsilon> 0$. If $A$ is
uniformly
elliptic and has a stationary ergodic law, then $u_\varepsilon$ is
known to
converge as $\varepsilon\to0$ to $u_{\mathsf{h}}$, the solution of
\[
(1-\nabla\cdot A_{\mathsf{h}}\nabla ) u_{\mathsf{h}}= f,
\]
where $A_{\mathsf{h}}$ is the (constant in space, deterministic) \emph
{homogenized
matrix}. This asymptotic result becomes more interesting if we can:
\begin{longlist}[(2)]
\item[(1)] devise (provably) efficient techniques to compute the homogenized matrix;
\item[(2)] estimate the error in the convergence of $u_\varepsilon$ to
$u_{\mathsf{h}}$.
\end{longlist}
Doing so requires to introduce some additional assumption on the type
of correlations displayed by the random coefficients; we assume from
now on that they have a finite range of dependence. These problems were
discussed in several works
\cite{yuri,capiof,boupia,boivin,vardecay,conspe,confah} (see also
\cite{cafsou,armsma} for nondivergence form operators), but optimal error bounds
were worked out only recently in \cite{glotto,glotto2,glmo,gno} for
(1), and in \cite{homog,gno-2scale} for (2) (in the discrete-space setting).

While controlling the size of the errors in homogenization is useful,
it would be better (and it is our aim) to describe precisely what the
errors look like when $\varepsilon$ is small. As an analogy, if the
convergence of $u_\varepsilon$ to $u_{\mathsf{h}}$ is a law of large
numbers, then we are
looking for a central limit theorem.

The present paper is a first step toward this goal. In a discrete-space
setting, it was proved in \cite{glotto} that stationary correctors
exist for $d \geqslant3$ (recall that we assume that the random coefficients
have a finite range of dependence). In this case, let us write $\phi
_\xi$ for the (stationary) corrector in the direction $\xi$ [see (\ref{defcorr})]. Under a minor smoothness assumption on the random
coefficients, we show that for large~$x$, the correlation
$ \langle\phi_\xi(0) \tarpas  \phi_\xi(x)  \rangle$
becomes very close to
%
\begin{equation}
\label{defK}
\mathscr{K}_\xi(x) := \int\nabla\mathcal{G}_{\mathsf{h}}(y)
\cdot \mathsf{Q}^{(\xi)}\, \nabla\mathcal{G}_{\mathsf{h}}(y-x) \tarpas\, \mathrm{d}y,
\end{equation}
where $\mathcal{G}_{\mathsf{h}}$ is the Green function of the
homogenized operator $-\nabla
\cdot A_{\mathsf{h}}\nabla$, and $\mathsf{Q}^{(\xi)}$ is a $d\times
d$ symmetric matrix that can be
expressed in terms of correctors, see (\ref{defQx}).

This result paves the way for the understanding of the full scaling
limit of $\varepsilon^{-({d-2})/{2}} \phi(\cdot/\varepsilon)$,
seen as a random
distribution. Indeed, the main missing ingredient is now to show that
for any bounded, smooth test function $f$, the properly rescaled random
variable $\sum_{x \in\mathbb{Z}^d} \phi(x) f(\varepsilon x)$
converges in law to a
Gaussian. This will be done in \cite{fluct}.

This result on the corrector suggests (via a formal two-scale
expansion) a scaling limit for $\varepsilon^{-d/2}
(u_\varepsilon(\cdot
)- \langle u_\varepsilon(\cdot)  \rangle )$ as
well. This will be addressed in \cite{fluct2}.

\subsection*{Related works}
We now give a brief overview of
related works. These can be divided into three groups.

First, the questions that we consider here in dimension $d \geqslant3$ have
been investigated in dimension $1$. One can benefit from this setting
to gain a better understanding of the effect of long-range correlations
of the coefficients \cite{balgar,gubal-1d}.

Second, similar questions have been explored for the homogenization of
operators other than those considered here. Typically, one considers a
deterministic operator perturbed by the addition of a rapidly
oscillating random potential \cite{fop,bal1,bal2,bal3,baljin,balgar2,gubal}. We refer to \cite{balgu-rev} for
a review.

Third, there is a deep connection between the corrector studied in the
present paper and so-called $\nabla\varphi$ interface models \cite{fun}.
At a heuristic level, one can think of the corrector as the
zero-temperature limit of such an interface model (with a
bond-dependent potential). The scaling limit of the interface model
with convex, homogeneous potential was shown to be the Gaussian free
field \cite{nadspe,gos,mill}. In view of this, one may expect (as was
suggested in \cite{berbis}, Conjecture~5) the correlations of the
corrector to be described by a Gaussian free field as well. However,
our results show that such is not the case in general. One way to see
this is to observe that the Fourier transform of $\mathscr{K}_\xi$ is
%
\begin{equation}
\label{efourier}
\frac{p \cdot\mathsf{Q}^{(\xi)}p}{(p \cdot A_{\mathsf{h}}p)^2} \qquad\bigl(p \in\mathbb{R}^d\bigr),
\end{equation}
while it should be of the form
%
\begin{equation}
\label{efourier-gff}
\frac{1}{p \cdot B p} \qquad\bigl(p \in\mathbb{R}^d\bigr)
\end{equation}
for some symmetric, positive definite matrix $B$, if the correlations
were those of a Gaussian free field. By considering coefficients with
small ellipticity ratio, one can produce examples where (\ref
{efourier}) cannot be reduced to (\ref{efourier-gff}).

The proof given in \cite{nadspe} that the interface model rescales to
the Gaussian free field (and the proof of the dynamical version of this
in \cite{gos}) uses a Helffer--Sj\"{o}strand representation of the
correlations. We will also use this representation here, but with an
important difference. In the case of the interface model, the
Helffer--Sj\"{o}strand representation readily enables to express the
correlations of the interface as the averaged Green function of some
operator, and the crux is then to show that this operator can be
homogenized. In our case, the representation has a less clear
interpretation. But it has to be so, since otherwise this would lead to
Gaussian-free-field correlations.

Recently, a very interesting and direct connection was put forward in
\cite{bisspo} between certain interface models with homogeneous but
possibly nonconvex potentials and the corrector considered here. The
authors obtained the scaling limit of interface models with such
potentials and zero tilt. They point out that the understanding of
models with nonzero tilt could be obtained from the understanding of
the scaling limit of the corrector. We refer to \cite{bisspo}, Section~6,
for more on this.

\subsection*{Organization of the paper}
The precise setting and
results of this paper are laid down in the next section. The
Helffer--Sj\"{o}strand representation of correlations is introduced in
Section~\ref{sHJ}. Section~\ref{sgreen} recalls several crucial
estimates on the corrector and the Green function. The goal of
Section~\ref{stwo-scale} is to justify, in a weak sense, the two-scale
expansion of the gradient of the Green function. The proof of the main
result is then completed in Section~\ref{sproof}.

\section{Precise setting and results}
\label{ssetting}

We consider the (nonoriented) graph $(\mathbb{Z}^d, \mathbb{B})$
with $d \geqslant3$, where
$\mathbb{B}$ is the set of nearest-neighbor edges. Let $(\mathbf
{e}_1,\ldots
,\mathbf{e}_d)$ be
the canonical basis of $\mathbb{Z}^d$. For every edge $e \in\mathbb
{B}$, there exists
a unique pair $(\underline{e},i) \in\mathbb{Z}^d \times\{1,\ldots,d\}$ such that $e$
links $\underline{e}$ to $\underline{e} + \mathbf{e}_i$. Given such
a pair, we
write $\overline{e}
= \underline{e} + \mathbf{e}_i$. We call $\underline{e}$ the \emph
{base point} of the
edge~$e$. For $f : \mathbb{Z}^d \to\mathbb{R}$, we let $\nabla f :
\mathbb{B}\to\mathbb{R}$ be the
gradient of $f$, defined by
\[
\nabla f(e) = f(\overline{e}) - f(\underline{e}).
\]
We write $\nabla^*$ for the formal adjoint of $\nabla$, that is, for
$F :
\mathbb{B}
\to\mathbb{R}$, $\nabla^*F : \mathbb{Z}^d \to\mathbb{R}$ is
defined via
\[
\bigl(\nabla^* F\bigr) (x) = \sum_{i = 1}^d
F\bigl((x-\mathbf{e}_i,x)\bigr) - F\bigl((x,x+\mathbf{e}_i)
\bigr).
\]
For such $F$, we define $AF(e) = a_e F(e)$, where $(a_e)$ are real
numbers taking values in a compact subset of $(0,+\infty)$. The
operator of interest is $\nabla^* A \nabla$.

While a standard assumption for our purpose would be that $(a_e)$ are
independent and identically distributed, the technicalities of the
proof will be reduced by assuming that they are also smooth in the
following sense. We give ourselves a family $(\zeta_e)_{e \in\mathbb
{B}}$ of
independent standard Gaussian random variables (we write $\mathbb{P}$
for the
law of this family on $\Omega= \mathbb{R}^{\mathbb{B}}$, and $
\langle\cdot  \rangle$ for the
associated expectation). The coefficients $(a_e)_{e \in\mathbb{B}}$
are then
defined by $a_e = a(\zeta_e)$, where $a : \mathbb{R}\to\mathbb{R}$
is a fixed twice
differentiable function with bounded first and second derivatives [and
taking values in a compact subset of $(0,+\infty)$].

Under these conditions, it is well known that there exists a constant
matrix $A_{\mathsf{h}}$ such that $\nabla^* A \nabla$ homogenizes
over large scales to
the continuous operator $-\nabla\cdot A_{\mathsf{h}}\nabla$.

Let $\xi$ be a fixed vector of $\mathbb{R}^d$. For $\mu> 0$, let
$\phi_{\xi,\mu}$ be
the unique stationary solution of
%
\begin{equation}
\label{defcorrmu}
\mu\phi_{\xi,\mu}+ \nabla^* A (\xi+ \nabla
\phi_{\xi,\mu}) = 0.
\end{equation}
It is proved in \cite{glotto} that (recall that we assume $d \geqslant3$)
$\phi_{\xi,\mu}$ converges in $L^2(\Omega)$ to the unique
stationary solution
$\phi_\xi$ of
%
\begin{equation}
\label{defcorr}
\nabla^* A (\xi+ \nabla\phi_\xi) = 0.
\end{equation}
The function $\phi_\xi$ is called the (stationary) \emph{corrector} in
the direction $\xi$. We use $\phi_i$ as shorthand for $\phi_{\mathbf
{e}_i}$. In
equations such as (\ref{defcorr}), $\xi$ is to be understood as the
function from $\mathbb{B}$ to $\mathbb{R}$ such that $\xi(e) = \xi
\cdot(\overline
{e}-\underline{e})$.

Let $\partial_e$ denote the weak derivative with respect to the random
variable $\zeta_e$, which we may call a \textit{vertical derivative}. The
formal adjoint of $\partial_e$ is
\[
\partial^*_e = -\partial_e + \zeta_e.
\]
We write $\partial f = (\partial_e f)_{e \in\mathbb{B}}$. For $F =
(F_e)_{e \in\mathbb{B}}$, we
write $\partial^* F = \sum_e \partial^*_e F_e$, and we let
%
\begin{equation}
\label{defL}
\mathscr{L}= \partial^*\, \partial.
\end{equation}
We write $|x|$ for the $L^2$-norm of $x \in\mathbb{Z}^d$. In order to keep
light notation, we let $|x|_* = |x|+ 2$ (so that, e.g., $\log
|x|_*$ is bounded away from $0$).

Here is our main result.

\begin{thm}[(Structure of correlations)]
\label{tstruct}
Recall that we assume $d \geqslant3$. Let $\mathcal{E}_0$ be the set
of edges
with base-point $0 \in\mathbb{Z}^d$, let $\mathcal{G}_{\mathsf{h}}:
\mathbb{R}^d \to\mathbb{R}$ be the Green
function of the (continuous-space) homogenized operator $- \nabla\cdot
A_{\mathsf{h}}\nabla$, let $\mathsf{Q}^{(\xi)}= (\mathsf{Q}^{(\xi
)}_{jk})_{1\leqslant j,k\leqslant d}$ be the matrix
defined by
%
\begin{eqnarray}
 \mathsf{Q}^{(\xi)}_{jk}
&=&\sum_{e \in\mathcal{E}_0} \bigl\langle(\mathbf{e}_j
+ \nabla\phi _j) (e) (\xi + \nabla\phi _\xi) (e) \tarpas
\,\partial_e a_e (\mathscr{L}+ 1)^{-1}
\nonumber
\\[-8pt]
\label{defQx}
\\[-8pt]
\nonumber
&&{}\times\partial_e a_e \tarpas (\mathbf{e}_k +
\nabla \phi_k) (e) (\xi+ \nabla\phi_\xi) (e) \bigr\rangle,\nonumber
\end{eqnarray}
and let $\mathscr{K}_\xi(x)$ be defined by (\ref{defK}). There
exists a constant $C
< \infty$ such that for every $x \in\mathbb{Z}^d \setminus\{0\}$,
%
\begin{equation}
\label{estruct}
\bigl\llvert \bigl\langle\phi_\xi(0) \tarpas
\phi_\xi(x) \bigr\rangle - \mathscr{K}_\xi(x) \bigr\rrvert
\leqslant C \tarpas \frac{\log
^2 |x|_*}{|x|^{d-1}}.
\end{equation}
\end{thm}

\begin{rem}
When $\mathbf{e}_j$ is interpreted as a function over $\mathbb{B}$ as in
(\ref{defQx}),
it is to be understood as $\xi$ is in (\ref{defcorr}), that is,
$\mathbf{e}
_j(e)$ is $1$ if the edge $e$ is parallel to the basis vector $\mathbf{e}_j$,
and is $0$ otherwise.
\end{rem}

\begin{rem}
The operator $\mathscr{L}$ is the infinitesimal generator of the
Ornstein--Uhlenbeck semigroup on $\mathbb{R}^{\mathbb{B}}$, and
$\mathbb{P}$ is a reversible
measure for the associated dynamics. For more general distributions of
coefficients, one may replace $\mathscr{L}$ by the infinitesimal
generator of
the Glauber dynamics, that is, to keep the definition~(\ref{defL}), but
with $\partial_e$ changed for
\[
\partial_e f= \mathbb{E}\bigl[f  |  (a_{e'})_{e' \neq e}
\bigr] - f
\]
[in which case $(\mathscr{L}+ 1)^{-1}$ must be replaced by $\mathscr
{L}^{-1}$ in (\ref{defQx})]. The setting we have chosen reduces the
amount of technicality
mostly by allowing us to use the chain rule for derivation.
\end{rem}

\begin{rem}
We learn from Proposition~\ref{pcontract} and Theorem~\ref
{tintegrability} that the tensor $\mathsf{Q}^{(\xi)}$ is well defined.
From the identity
%
\begin{eqnarray}
&& \xi'\cdot\mathsf{Q}^{(\xi)}
\xi'
\nonumber
\\
\label{F1}
&&\qquad=\sum_{e \in\mathcal{E}_0} \bigl\langle\bigl(
\xi' + \nabla\phi_{\xi
'}\bigr) (e) (\xi+ \nabla\phi
_\xi) (e) \tarpas \\
&&\qquad\quad{}\times\partial_e a_e (
\mathscr{L}+ 1)^{-1} \,\partial_e a_e \tarpas
\bigl(\xi' + \nabla\phi_{\xi'}\bigr) (e) (\xi+ \nabla
\phi_\xi) (e) \bigr\rangle,\nonumber
\end{eqnarray}
which follows from the linearity of $\phi_\xi$ in $\xi$,
we learn that $\mathsf{Q}^{(\xi)}$ is positive semi-definite.
In particular, the Fourier transform of $\mathscr{K}_\xi$ is nonnegative.

Moreover, $\mathsf{Q}^{(\xi)}$ is nondegenerate as soon as the
derivative of the
function $a : \mathbb{R}\to\mathbb{R}$ is everywhere positive.
Indeed, if the expression (\ref{F1}) vanishes for $\xi'=\xi$, the
strict positivity of the operator implies that for any $e\in\mathcal{E}_0$,
$\partial_e a_e(\xi+ \nabla\phi_\xi)^2(e)$ and thus $(\xi+ \nabla
\phi_\xi
)^2(e)$ vanishes almost surely. This in turn implies that
\[
\sum_{e \in\mathcal{E}_0} \bigl\langle(\xi+ \nabla
\phi_\xi ) (e)a_e(\xi+ \nabla \phi _\xi) (e)
\bigr\rangle=\xi\cdot A_h\xi
\]
vanishes. By the nondegeneracy of the homogenized tensor $A_h$, this
yields as desired $\xi=0$. The same argument also implies that
the null space of $\mathsf{Q}^{(\xi)}$ is contained in the hyperplane
orthogonal to
$A_h\xi$.
\end{rem}
%
\begin{rem}
There is no simple relation between the quartic form defined by
$\mathsf{Q}^{(\xi)}
_{jk}$ and the quadratic form $A_h$,
besides that $\xi'\mathsf{Q}^{(\xi)}\xi'$ is bounded from below by
$(\xi'\cdot
A_h\xi
)^2$ up to a multiplicative constant.
As was noted in the \hyperref[sintro]{Introduction}, $\mathscr{K}_\xi$ is not the Green
function of a
second-order operator in general. While its Fourier transform has
the right sign and homogeneity, it is not the inverse of a quadratic form.
\end{rem}
%
\begin{rem}
By polarization of the quartic form $\mathsf{Q}^{(\xi)}_{jk}$ in the
$\xi$-variables,
one also obtains a result for covariances $ \langle\phi_\xi(0)
\tarpas  \phi _{\xi'}(x)  \rangle$ with $\xi'\neq\xi$.
\end{rem}


\begin{rem}
We expect that at least if the environment is sufficiently mixing (as,
e.g., when its correlations are of finite range), then there
exists a matrix~$\mathsf{Q}^{(\xi)}$ [whose explicit expression may
differ from that
given in (\ref{defQx})] such that the large-scale correlations of the
correctors are described by (\ref{defK}) and (\ref{estruct}).
\end{rem}

%
%
\section{Hellfer--Sj\"ostrand representation}
\label{sHJ}
%
\begin{prop}[(Helffer--Sj\"{o}strand representation of correlations,
\cite{helsjo,sjo,nadspe})]\label{pHJ}
Let $f,g : \Omega\to\mathbb{R}$ be centered square-integrable
functions such
that for every $e \in\mathbb{B}$, $\partial_e f$ and $\partial_e g$
are in
$L^2(\Omega
)$. We have
\[
\langle f \tarpas g \rangle = \sum_{e \in\mathbb{B}} \bigl
\langle\partial_e f \tarpas (\mathscr{L}+1)^{-1} \tarpas
\,\partial_e g \bigr\rangle.
\]
\end{prop}
\begin{pf}
The claim is similar to (and simpler than) that obtained in
\cite{nadspe}, Section~2.1. We recall the proof briefly for the reader's
convenience. By density, we can restrict our attention to functions $f,
g$ that depend only on a finite number of $(\zeta_e)_{e \in\mathbb
{B}}$, and
also by density, we may assume $f$ and $g$ to be smooth functions.
Note that the commutator $[\partial_e, \partial_{e'}^*]$ satisfies
%
\begin{equation}
\label{commutator}
\bigl[\partial_e, \partial_{e'}^*\bigr] =
\mathbf{1}_{e = e'}.
\end{equation}

Let us first assume that there exists a function $u \in L^2(\Omega)$
such that
$g = \mathscr{L}u$. Writing $G = \partial u$, we observe that
\begin{eqnarray*}
\partial_e g & = & \partial_e\, \partial^* G
\\
& = & \sum_{e'} \partial_e\,
\partial_{e'}^* G_{e'}
\\
& = & \sum_{e'} \bigl(\bigl[\partial_e,
\partial_{e'}^*\bigr] + \partial _{e'}^* \,\partial
_e \bigr) G_{e'}
\\
& = & G_e + \sum_{e'}
\partial_{e'}^* \,\partial_{e'} G_e,
\end{eqnarray*}
where we used (\ref{commutator}) and the fact that $\partial_{e'} G_e =
\partial_{e'} \,\partial_e u = \partial_e G_{e'}$ in the last step.
Recalling the definition of $\mathscr{L}$ in (\ref{defL}), we arrive at
\[
\partial_e g = (\mathscr{L}+ 1) G_e = (\mathscr{L}+ 1)
\,\partial_e u.
\]
In particular, $\partial_e u \in L^2(\Omega)$ and
\[
\langle f \tarpas g \rangle = \langle f \tarpas \mathscr{L}u \rangle = \sum
_e \langle\partial_e f \tarpas
\,\partial_e u \rangle =\sum_e \bigl
\langle\partial_e f \tarpas (\mathscr {L}+1)^{-1}
\,\partial_e g \bigr\rangle.
\]
In order to conclude, it suffices to check that 
the range of the operator $\mathscr{L}$ is dense in the set of centered
square-integrable functions. If $g \in L^2(\Omega)$ is smooth, depends
on a finite number of $(\zeta_e)_{e \in\mathbb{B}}$ and is in the orthogonal
complement of $\operatorname{Ran}(\mathscr{L})$, then
\[
\langle g \tarpas \mathscr{L}g \rangle = 0 = \sum_e
\bigl\langle|\,\partial_e g|^2 \bigr\rangle,
\]
so $g$ is constant. It follows that the orthogonal complement of
$\textrm{Ran}(\mathscr{L})$ is the set of constant functions, and
this completes
the proof.
\end{pf}
The following additional information on $(\mathscr{L}+1)^{-1}$ will
turn out to
be useful.
%
\begin{prop}[(Contraction of $L^p$)]
\label{pcontract}
For every $p \geqslant2$, the operator $(\mathscr{L}+1)^{-1}$ is a
contraction from
$L^p(\Omega)$ to itself.
\end{prop}
\begin{pf}
Let $\Lambda$ be a finite subset of $\mathbb{B}$, and let $\mathcal
{F}_\Lambda$ the
set of real functions of $(\zeta_e)_{e \in\Lambda}$. We define
$H^1_\Lambda$ as the completion of the set of smooth functions in
$\mathcal
{F}_\Lambda$ for the scalar product
\[
(u,v)_{H^1_\Lambda} = \langle u \tarpas v \rangle + \sum
_{e
\in\Lambda} \langle\partial _e u \tarpas
\,\partial_e v \rangle.
\]
For every $f \in H^1_\Lambda$, there exists a unique $u \in
H^1_\Lambda
$ such that
%
\begin{equation}
\label{weak}
\forall v \in H^1_\Lambda,\qquad \tarpas
(u,v)_{H^1_\Lambda} = \langle f \tarpas v \rangle,
\end{equation}
and this is nothing but the weak formulation of the equation $(\mathscr
{L}+1)u =
f$. For every $\varepsilon> 0$, let
$
\psi_\varepsilon(x) = \varepsilon^{-1} \arctan(\varepsilon x)
$
be a ``nice'' (in particular, bounded) approximation of the function $x
\mapsto x$.
One can check that if $v \in H^1_\Lambda$, then $\psi_\varepsilon(v
\tarpas  |v|^{p-2}) \in H^1_\Lambda$. Hence, for $u \in H^1_\Lambda$ satisfying
(\ref{weak}),
\[
\bigl(u,\psi_\varepsilon\bigl(u \tarpas |u|^{p-2}\bigr)
\bigr)_{H^1_\Lambda} = \bigl\langle f \tarpas \psi_\varepsilon \bigl(u
\tarpas |u|^{p-2}\bigr) \bigr\rangle,
\]
and we recall that
\[
\bigl(u,\psi_\varepsilon\bigl(u \tarpas |u|^{p-2}\bigr)
\bigr)_{H^1_\Lambda} = \bigl\langle u \tarpas \psi_\varepsilon \bigl(u
\tarpas |u|^{p-2}\bigr) \bigr\rangle + \sum_{e \in
\Lambda}
\bigl\langle\partial_e u \tarpas \,\partial_e \psi
_\varepsilon \bigl(u \tarpas |u|^{p-2}\bigr) \bigr\rangle.
\]
Since $u \mapsto\psi_\varepsilon(u \tarpas  |u|^{p-2})$ is an increasing function,
it follows that for every $e$,
\[
\bigl\langle\partial_e u \tarpas \,\partial_e
\psi_\varepsilon\bigl(u \tarpas |u|^{p-2}\bigr) \bigr\rangle
\geqslant0,
\]
and thus
\[
\bigl\langle u \tarpas \psi_\varepsilon\bigl(u \tarpas |u|^{p-2}
\bigr) \bigr\rangle \leqslant \bigl\langle f \tarpas \psi_\varepsilon\bigl(u
\tarpas |u|^{p-2}\bigr) \bigr\rangle.
\]
By the monotone convergence theorem, the left-hand side converges to
$ \langle|u|^p  \rangle= \|u\|_p^p$ as $\varepsilon$
tends to $0$. The right-hand side is bounded by
\[
\|f\|_p \tarpas \bigl\langle\bigl|\psi_\varepsilon\bigl(u \tarpas
|u|^{p-2}\bigr)\bigr|^{p/({p-1})} \bigr\rangle^{1-1/p}\leqslant\| f
\| _p\|u\|_p^{p-1},
\]
where we have used $|\psi_\varepsilon(x)|\leqslant|x|$. We have thus shown
\[
\|u\|_p^p \leqslant\|f\|_p \tarpas \|u
\|_p^{p-1},
\]
that is, $\|u\|_p \leqslant\|f\|_p$, and this implies the theorem.
\end{pf}
Using the fact that $(\mathscr{L}+1)^{-1}$ is a contraction on
$L^2(\Omega)$, we
deduce the following covariance estimate, which parallels those
appearing in \cite{nadspe,nadspe-unpub} (Brascamp--Lieb inequality),
\cite{gno}, Definition~1  and \cite{glotto2}, Lemma~3.
%
\begin{cor}[(Covariance estimate)]
\label{ccovar}
For $f$ and $g$ as in Proposition~\ref{pHJ},
\[
\bigl\llvert \langle f \tarpas g \rangle \bigr\rrvert \leqslant\sum
_{e
\in\mathbb{B}} \bigl\langle(\partial_e f)^2
\bigr\rangle^{1/2} \bigl\langle(\partial_e g)^2
\bigr\rangle^{1/2}.
\]
\end{cor}
%
%
\section{Estimates on the corrector and the Green function}
\label{sgreen}
The aim of this section is to gather several known estimates on the
Green function and on the corrector.

\begin{thm}[(Existence and integrability of the corrector \cite
{glotto})]\label{tintegrability}\label{tcorrector}
Recall that we assume $d \geqslant3$. For every $\mu> 0$, there
exists a
unique stationary solution $\phi_{\xi,\mu}$ to equation (\ref{defcorrmu}).
Moreover, for every $p \geqslant1$, $ \langle|\phi_{\xi,\mu
}(0)|^p  \rangle$ and $ \langle|\nabla\phi_{\xi,\mu
}(e)|^p  \rangle$ ($e \in\mathbb{B}$) are uniformly bounded in
$\mu> 0$.
The limit
\[
\phi_\xi= \lim_{\mu\to0} \phi_{\xi,\mu}
\]
is well defined in $L^p(\Omega)$ and is the unique stationary solution
to (\ref{defcorr}).
\end{thm}
A direct consequence of this result is:
%
\begin{cor}[(Almost-sure control of the corrector)]\label{ccorrector-growth}
Let $B_n = \{
-n,\break \ldots,  n\}^d$ and let $\mathbb{B}_n$ be the set of edges whose base-point is in
$B_n$. For every $\beta> 0$, almost surely,
\[
\lim_{n \to+\infty} \tarpas n^{-\beta} \max_{x \in B_n}
\bigl|\phi_\xi(x)\bigr| = 0
\]
and
\[
\lim_{n \to+\infty} \tarpas n^{-\beta} \max_{e \in\mathbb{B}_n}
\bigl|\nabla\phi_\xi (e)\bigr| = 0.
\]
\end{cor}
\begin{pf}
Let $p \geqslant1$. By Chebyshev's inequality,
\[
\mathbb{P}\bigl[\bigl|\phi_\xi(0)\bigr| \geqslant x\bigr] \leqslant
\frac{\mathbb
{E}[|\phi_\xi(0)|^p]}{x^p} \qquad(x > 0),
\]
so for any $\varepsilon> 0$, by a union bound,
\[
\mathbb{P} \Bigl[n^{-\beta} \max_{x \in B_n} \bigl|
\phi_\xi(x)\bigr| \geqslant\varepsilon \Bigr] \leqslant|B_n|
\frac{\mathbb{E}[|\phi_\xi(0)|^p]}{(\varepsilon n^\beta)^p}.
\]
The first part of the corollary follows by taking $p$ large enough and
applying the Borel--Cantelli lemma. The second part is obtained in the
same way.
\end{pf}

We write $G(x,y)$ for the Green function between points $x$ and $y$ in
$\mathbb{Z}^d$, that is, $G(x,y) =  (\nabla^* A \nabla
)^{-1}(x,y)$ [the
dependence on $(a_e)_{e \in\mathbb{B}}$ is kept implicit in the
notation]. For
$\mu> 0$, we also let $G_\mu(x,y) =  (\mu+ \nabla^* A \nabla
 )^{-1}(x,y)$.

Regularity theory ensures the following decay properties of the Green
function (see, e.g., \cite{dl-diff}, Proposition~3.6, for a
proof adapted to our context).
%
\begin{thm}[(Pointwise estimates on the Green function)]\label{tpointGreen}
There exist $C<\infty$, $c > 0$ and $\alpha> 0$ such that for every
$\mu\in[0,1/2]$ and $\zeta\in\Omega$,
%
\begin{eqnarray}
\label{pointGreen}
G_\mu(0,x)  &\leqslant & \frac{C}{|x|_*^{d-2}} \tarpas
e^{-{c} \sqrt{\mu} |x|} \qquad \bigl(x \in\mathbb{Z}^d\bigr),
\\
\label{pointgradGreen}
\bigl\llvert \nabla G_\mu(0,e)\bigr\rrvert  &\leqslant &
\frac{C}{|\underline
{e}|_*^{d-2+\alpha}} \tarpas e^{-{c}
\sqrt{\mu} |\underline{e}|} \qquad(e \in\mathbb{B}).
\end{eqnarray}
\end{thm}
It was recently shown in \cite{marott} that, after averaging over the
environment, the rates of decay of the gradient and mixed second
gradient of the Green function behave as in the homogeneous case (see
also \cite{dl-diff}, Remark~11.2,  for the fact that the estimates hold
uniformly over $\mu$).
%
\begin{thm}[(Annealed estimates on the gradients of the Green function
\cite{marott})]\label{tmarott}
For every $1\leqslant p <\infty$, there exists $C<\infty$ such that
for every
$\mu\in[0,1/2]$ and every $e,e' \in\mathbb{B}$,
\begin{eqnarray*}
\bigl\langle\bigl\llvert \nabla G_\mu(0,e) \bigr\rrvert
^p \bigr\rangle ^{1/p} & \leqslant & \frac{C}{|\underline{e}|_*^{d-1}},
\\
\bigl\langle\bigl\llvert \nabla\nabla G_\mu\bigl(e,e'
\bigr) \bigr\rrvert ^p \bigr\rangle^{1/p}  &\leqslant &
\frac{C}{|\underline{e}'-\underline{e}|_*^{d}}.
\end{eqnarray*}
\end{thm}
%
\begin{rem}
Notice that $\nabla G(x,e)$ (for $x \in\mathbb{Z}^d$ and $e \in
\mathbb{B}$) denotes the
gradient of $G(x,\cdot)$ evaluated at the edge $e$. Similarly, $\nabla
\nabla
G(e,e')$ denotes the gradient of $\nabla G(\cdot,e')$ evaluated at the
edge $e$.
\end{rem}
We conclude this section by recalling useful computations of vertical
derivatives. The following two propositions are borrowed from \cite{glotto}, Lemmas~2.4 and 2.5.
%
\begin{prop}[(Derivatives of the corrector \cite{glotto})]
\label{pderiv-corr}
For every $\mu> 0$, $x \in\mathbb{Z}^d$ and $e \in\mathbb{B}$, the
approximate
corrector $\phi_{\xi,\mu}(x)$ is differentiable with respect to
$\zeta_e$ and
\[
\partial_e \phi_{\xi,\mu}(x) = - \partial_e
a_e \nabla G_\mu(x,e) (\xi+ \nabla\phi_{\xi,\mu})
(e).
\]
\end{prop}
%
\begin{rem}
Recalling that we assume $a_e$ to be of the form $a(\zeta_e)$ with $a$
differentiable, we can rewrite $\partial_e a_e$ as $a'(\zeta_e)$.
\end{rem}
%
\begin{rem}
\label{rphixmphimu}
Contrary to $\phi_{\xi,\mu}$, the corrector $\phi_\xi$ is not
well defined for
every value of $(\zeta_e) \in\Omega$, but only on a set of full
probability measure. In order to prove a statement similar to
Proposition~\ref{pderiv-corr} for $\phi_\xi$ instead of $\phi_{\xi
,\mu}
$, it
is thus necessary to show first that $\phi_\xi$ is defined on a subset
of $\Omega$ large enough that speaking of $\partial_e \phi_\xi$ be
meaningful. We will however not show this here, since for our purpose,
it is always possible to bypass this problem by approximating $\phi
_\xi
$ by $\phi_{\xi,\mu}$, computing the derivatives, and then passing
to the
limit $\mu\to0$.
\end{rem}
%
\begin{prop}[(Derivatives of the Green function \cite{glotto})]\label{pderiv-green}
For every $\mu\geqslant0$, $x,y \in\mathbb{Z}^d$ and $e \in\mathbb
{B}$, the Green
function $G_\mu(x,y)$ is differentiable with respect to $\zeta_e$ and
\[
\partial_e G_\mu(x,y) = - \partial_e
a_e \nabla G_\mu(x,e) \nabla G_\mu(y,e).
\]
\end{prop}
These two propositions can be proved by differentiating the defining
equation of, respectively, the corrector and the Green function, namely
\begin{eqnarray*}
\mu\phi_{\xi,\mu}+ \nabla^* A (\xi+ \nabla\phi_{\xi,\mu}) &=& 0,
\\
\bigl(\mu+ \nabla^* A \nabla\bigr) G_\mu(x,\cdot) &=&
\mathbf{1}_x.
\end{eqnarray*}
We refer to \cite{glotto} for details.

\section{Two-scale expansion of the Green function}
\label{stwo-scale}
Note that since we assume the coefficients to be independent and
identically distributed, the law of the coefficients is invariant under
the rotations that preserve the lattice, and $A_{\mathsf{h}}$ is thus
a multiple
of the identity, say $A_{\mathsf{h}}= a_{\mathsf{h}}\mathrm{Id}$. We
define the discrete
homogenized Green function $G_{\mathsf{h}}$ as the unique bounded
solution of the equation
\[
\nabla^* A_{\mathsf{h}}\nabla G_{\mathsf{h}}= \mathbf{1}_0,
\]
where $A_{\mathsf{h}}$ in the formula above acts as the multiplication by
$a_{\mathsf{h}}$ on
every edge.
For $f : \mathbb{Z}^d \to\mathbb{R}$ and $x \in\mathbb{Z}^d$, we
write $\nabla_j f(x)$ to denote
$f(x+\mathbf{e}_j) - f(x)$. If instead we take $e \in\mathbb{B}$, we
understand $\nabla_j
f(e)$ to mean $\nabla_j f(\underline{e})$, that is, the gradient of
$f$ along the
edge parallel to the vector $\mathbf{e}_j$ having the same base-point
as $e$.

The goal of this section is to prove the following quantitative
two-scale expansion of the gradient of the Green function.
%
\begin{thm}[(Quantitative two-scale expansion of the Green function)]
\label{ttwo-scale}
For every $p > 2$, there exists $C<\infty$ such that the following
holds. If $g : \Omega\to\mathbb{R}$ is in $L^p(\Omega)$ and is
differentiable
with respect to $\zeta_b$ with $\partial_b g \in L^p(\Omega)$ for
every $b
\in\mathbb{B}$, then for every $e \in\mathbb{B}$,
%
\begin{eqnarray}
&& \Biggl\llvert \bigl\langle g \tarpas \nabla G(0,e) \bigr\rangle -
\sum_{j =
1}^d \nabla_j
G_{\mathsf{h}}(e) \bigl\langle g \tarpas (\mathbf{e} _j + \nabla
\phi_j) (e) \bigr\rangle \Biggr\rrvert
\nonumber
\\[-8pt]
\label{etwo-scale}
\\[-8pt]
\nonumber
&&\qquad\leqslant C \biggl( \|g\|_p \frac{\log|\underline{e}|_*}{|\underline
{e}|_*^d} + \mathop{\sum_{y \in\mathbb{Z}^d}}_{b \in\mathbb{B}} \llVert \partial_b
g \rrVert _p \frac{1}{|\underline
{e}-y|_*^{d-1} \tarpas  |\underline{b}-y|_*^{d} \tarpas  |y|_*^d} \biggr).
\end{eqnarray}
\end{thm}
%
\begin{rem}
Applying Theorem~\ref{ttwo-scale} with $g = 1$, we obtain that
\[
\bigl\llvert \bigl\langle\nabla G(0,e) \bigr\rangle - \nabla G_{\mathsf
{h}}(e)
\bigr\rrvert \leqslant C \tarpas \frac{\log|\underline{e}|_*}{|\underline{e}|_*^d}.
\]
\end{rem}
%
\begin{rem}
\label{rtranslation}
By translation, under the assumptions of Theorem~\ref{ttwo-scale}, we
also have
\begin{eqnarray*}
&& \Biggl\llvert \bigl\langle g \tarpas \nabla G(x,e) \bigr\rangle - \sum
_{j =
1}^d \nabla_j
G_{\mathsf{h}}(e-x) \bigl\langle g \tarpas (\mathbf{e}_j + \nabla
\phi_j) (e) \bigr\rangle \Biggr\rrvert
\\
&&\qquad\leqslant C \biggl( \|g\|_p \frac{\log|\underline
{e}-x|_*}{|\underline{e}-x|_*^d} + \mathop{
\sum_{y \in\mathbb{Z}^d }}_{ b \in\mathbb{B}} \llVert
\partial_b g \rrVert _p \frac
{1}{|\underline{e}-y|_*^{d-1} \tarpas  |\underline{b}-y|_*^{d}  |y-x|_*^d} \biggr),
\end{eqnarray*}
where $e-x$ denotes the translation of the edge $e$ by the vector $-x$.
\end{rem}
We define $z : \mathbb{Z}^d \to\mathbb{R}$ by
%
\begin{equation}
\label{defz}
z(x) = G(0,x) - G_{\mathsf{h}}(x) - \sum
_{j = 1}^d \phi_j(x) \tarpas
\nabla_j G_{\mathsf{h}}(x).
\end{equation}
%
\begin{prop}[(Equation for $z$ \cite{papvar,gno-2scale})]\label{pz-eq}
Let $A_i(x)$ stand for $a_{x,x+\mathbf{e}_i}$. Write $\nabla^2
G_{\mathsf{h}}$ for
the matrix
with entries $\nabla^*_i \nabla_j G_{\mathsf{h}}$ ($1 \leqslant i,j
\leqslant d$). Let $R$ be the
matrix with entries $(R_{ij})$ satisfying
\[
(R - A_{\mathsf{h}})_{ij} = - \bigl[A_i\bigl(
\mathbf{1}_j^i + \nabla_i \phi_j
\bigr) \bigr](\cdot- \mathbf{e}_i) \qquad (1\leqslant i,j\leqslant d),
\]
where $\mathbf{1}^i_{j} = \mathbf{1}_{i = j}$. For $e \in\mathbb
{B}$ in the direction of
$\mathbf{e}
_i$, let
\[
h(e) = - \Biggl(A \sum_{j = 1}^d
\phi_j(\cdot+ \mathbf{e}_i) \nabla \nabla _j
G_{\mathsf{h}} \Biggr) (e),
\]
and denote the $\mathbb{R}^{d\times d}$-scalar product of two matrices
$M$ and
$N$ by $M:N$ (i.e., the sum of all terms after entry-wise product).
We have
%
\begin{equation}
\label{ez-eq}
\nabla^* A \nabla z = R : \nabla^2 G_{\mathsf{h}}+
\nabla^* h.
\end{equation}
\end{prop}
%
\begin{rem}
\label{rcentredR}
The crucial feature of the right-hand side of (\ref{ez-eq}) is that it
involves only the second derivatives of $G_{\mathsf{h}}$ (this is
precisely what
one aims for when defining~$z$). Another aspect that will turn out to
be important for our purpose is that $ \langle R(x)
\rangle = 0$. This follows
from the fact [see, e.g., \cite{kuen}, (3.17)] that the $(i,j)$th
entry of the homogenized matrix $A_{\mathsf{h}}$ is equal to
\[
\bigl\langle A_i\bigl(\mathbf{1}_j^i +
\nabla_i \phi_j\bigr) \bigr\rangle = \bigl\langle
\mathbf{e}_i \cdot A(\mathbf{e}_j + \nabla
\phi_j) \bigr\rangle.
\]
\end{rem}
\begin{pf*}{Proof of Proposition~\protect\ref{pz-eq}}
We follow the line of argument given in the first step of the proof of
\cite{gno-2scale}, Theorem~1
(itself inspired by the first proof of \cite{papvar}, Theorem~3). For
$f : \mathbb{Z}^d \to\mathbb{R}$, we write $\nabla^*_i f(x) =
f(x-\mathbf{e}_i) - f(x)$.
To begin with, we observe that the following discrete Leibniz rules
hold, for $f,g : \mathbb{Z}^d \to\mathbb{R}$:
\begin{eqnarray*}
\nabla_i(fg) &=& (\nabla_i f) g + f(\cdot+
\mathbf{e}_i) \nabla_i g,
\\
\nabla^*_i(fg) &=& \bigl(\nabla^*_i f\bigr) g + f(\cdot-
\mathbf{e}_i) \nabla ^*_i g.
\end{eqnarray*}

Recall that by definition,
\[
\nabla^* A \nabla G(0,\cdot) = \mathbf{1}_0 = \nabla^*
A_{\mathsf
{h}}\nabla G_{\mathsf{h}},
\]
and thus,
\[
\nabla^* A \nabla\bigl(G(0,\cdot) - G_{\mathsf{h}}\bigr) = \nabla
^*(A_{\mathsf{h}}- A) \nabla G_{\mathsf{h}}.
\]
%
%
Writing $A_{\mathsf{h},i}$ for the $i$th diagonal coefficient of the (diagonal)
matrix $A_{\mathsf{h}}$, we can express the right-hand side above as
\[
\sum_{i = 1}^d \nabla^*_i(A_{\mathsf{h},i}-
A_i) \nabla_i G_{\mathsf{h}}.
\]
We now need to compute
%
\begin{equation}
\label{z-eq1}
\nabla^* A \nabla ( \phi_j\tarpas
\nabla_j G_{\mathsf{h}} ).
\end{equation}
By the Leibniz rule,
\[
\nabla_i ( \phi_j\tarpas \nabla_j
G_{\mathsf{h}} ) = (\nabla _i\phi_j)
\nabla_j G_{\mathsf{h}}+ \phi _j(\cdot+
\mathbf{e}_i) \nabla_i \nabla_j
G_{\mathsf{h}}.
\]
Hence, the term in (\ref{z-eq1}) is equal to
\[
\sum_{i = 1}^d \nabla^*_i
\bigl[ A_i \bigl(\nabla_i\phi_j \tarpas
\nabla _j G_{\mathsf{h}}+ \phi _j(\cdot+
\mathbf{e}_i) \nabla_i \nabla_j
G_{\mathsf{h}} \bigr) \bigr].
\]
We can thus rewrite $\nabla^* A \nabla z$ as
\begin{eqnarray*}
&&\sum_{i = 1}^d \Biggl\{
\nabla^*_i(A_{\mathsf{h},i}- A_i) \nabla_i
G_{\mathsf{h}}- \sum_{j = 1}^d \nabla
^*_i \bigl[ A_i \bigl(\nabla_i
\phi_j \nabla_j G_{\mathsf{h}}+ \phi _j(
\cdot+ \mathbf {e}_i) \nabla_i \nabla _j
G_{\mathsf{h}} \bigr) \bigr] \Biggr\}
\\
&&\qquad= \sum_{i=1}^d \Biggl
\{A_{\mathsf{h},i}\nabla^*_i \nabla_i G_{\mathsf{h}}\\
&&\qquad\quad{}-
\sum_{j = 1}^d \nabla ^*_i \bigl[
A_i \bigl(\bigl(\mathbf{1}^i_{j} +
\nabla_i\phi_j\bigr) \nabla_j
G_{\mathsf
{h}}+ \phi_j(\cdot+ \mathbf{e}_i) \nabla
_i \nabla_j G_{\mathsf{h}} \bigr) \bigr] \Biggr\},
\end{eqnarray*}
where we used the fact that $A_{\mathsf{h}}$ is constant. By the
definition of the
corrector, we have
\[
\sum_{i=1}^d \nabla^*_i
A_i\bigl(\mathbf{1}^i_{j} +
\nabla_i\phi_j\bigr) = \nabla^* A (\mathbf{e}_j
+ \nabla \phi_j) = 0,
\]
so by the Leibniz rule,
\[
\sum_{i,j=1}^d \nabla^*_i
\bigl[A_i\bigl(\mathbf{1}^i_{j} + \nabla
_i\phi_j\bigr) \nabla_j G_{\mathsf{h}}
\bigr] = \sum_{i,j=1}^d
\bigl[A_i\bigl(\mathbf{1}^i_{j} +
\nabla_i\phi_j\bigr) \bigr](\cdot- \mathbf
{e}_i) \tarpas \nabla ^*_i \nabla_j
G_{\mathsf{h}},
\]
and the conclusion follows.
\end{pf*}
As a consequence, we get the following representation for $z$.
%
\begin{prop}[(Representation for $z$)]
\label{prep-z}
For every $x \in\mathbb{Z}^d$,
%
\begin{equation}
\label{erep-z}
z(x) = \sum_{y \in\mathbb{Z}^d} G(x,y) \bigl( R :
\nabla ^2G_{\mathsf{h}} \bigr) (y) + \sum
_{b \in
\mathbb{B}} \nabla G(x,b) h(b).
\end{equation}
\end{prop}
\begin{pf}
Let $\tilde{z}(x)$ denote the right-hand side of (\ref{erep-z}),
which is
well defined by Corollary~\ref{ccorrector-growth}. Letting $\overline
{z} = z
- \tilde{z}$, one can check thanks to Proposition~\ref{pz-eq} that
$\nabla^*
A \nabla\overline{z} = 0$. In particular,
\[
\sum_{x \in B_n} \overline{z}(x) \nabla^* A \nabla
\overline{z}(x) = 0.
\]
This sum differs from
\[
\sum_{e \in\mathbb{B}_n} \nabla\overline{z}(e) \cdot A \nabla
\overline{z}(e)
\]
by no more than a constant times
%
\begin{equation}
\label{erepz1}
\sum_{e \in\mathbb{B}_{n+1} \setminus\mathbb{B}_{n}} \bigl(\bigl|\overline{z}(
\underline{e})\bigr| + \bigl|\overline {z}(\overline{e})\bigr| \bigr) \bigl|\nabla\overline{z}(e)\bigr|.
\end{equation}
This sum tends to $0$ as $n$ tends to infinity. To see this, we come
back to the definitions of $z$ and $\tilde{z}$, given, respectively, in
(\ref{defz}) and in the right-hand side of (\ref{erep-z}). Using
Corollary~\ref{ccorrector-growth}, Theorem~\ref{tpointGreen} and Proposition~\ref
{psums} of the   \hyperref[sappend]{Appendix}, we obtain that for every $\beta> 0$, almost surely,
\begin{eqnarray*}
\bigl|z(x)\bigr|  &=&  o \biggl(\frac{1}{|x|^{d-2-\beta}} \biggr) \qquad \bigl(|x| \to \infty\bigr),
\\
\bigl|\nabla{z}(e)\bigr| &=& o \biggl(\frac{1}{|\underline{e}|^{d-2+\alpha-\beta
}} \biggr)\qquad \bigl(|\underline{e}| \to
\infty\bigr)
\end{eqnarray*}
(where $\alpha$ comes from Theorem~\ref{tpointGreen}), and the same
relations hold for $z$ replaced by~$\tilde{z}$, and thus also for $z$
replaced by $\overline{z}$. Since $d \geqslant3$, we can take $\beta
> 0$
sufficiently small to ensure that $2(d-2)+\alpha-2\beta> d-1$, and we
obtain that the sum in (\ref{erepz1}) tends to $0$ as $n$ tends to infinity.

To sum up, we obtained that
\[
\lim_{n \to+\infty} \sum_{e \in\mathbb{B}_n} \nabla
\overline {z}(e) \cdot A \nabla \overline {z}(e) = 0.
\]
Since $A$ is positive definite, we conclude that $\overline{z}$ is a
constant. Now, both $z$ and $\tilde{z}$ tend to $0$ at infinity, so in
fact $\overline{z} = 0$, and this completes the proof.
\end{pf}
\begin{pf*}{Proof of Theorem~\ref{ttwo-scale}}
Let us first see that it suffices to show that
%
\begin{equation}
\label{ez-two}
\quad\bigl\llvert \bigl\langle g \tarpas \nabla z(e) \bigr\rangle
\bigr\rrvert \leqslant C \biggl( \|g\|_p \frac{\log|\underline{e}|_*}{|\underline
{e}|_*^d} + \mathop{
\mathop{\sum_{y \in\mathbb{Z}^d}}}_{ b \in\mathbb{B}} \llVert
\partial_b g \rrVert _p \frac{1}{|\underline
{e}-y|_*^{d-1} \tarpas  |y-\underline{b}|_*^{d} \tarpas  |y|_*^d} \biggr).
\end{equation}
Note that, by the Leibniz rule,
\[
\nabla_i z(x) = \nabla_i G(0,x) - \nabla_i
G_{\mathsf{h}}(x) - \sum_{j = 1}^d \bigl[
\nabla_i \phi_j(x) \tarpas \nabla_j
G_{\mathsf{h}}(x) + \phi_j(x+\mathbf{e}_i) \tarpas
\nabla_i \nabla_j G_{\mathsf{h}}(x) \bigr].
\]
In order to prove that (\ref{ez-two}) implies (\ref{etwo-scale}), it
is thus sufficient to show that
%
\begin{equation}
\label{z-3}
\bigl\llvert \bigl\langle g \tarpas \phi_j(x+
\mathbf{e}_i) \tarpas \nabla_i \nabla _j
G_{\mathsf{h}}(x) \bigr\rangle \bigr\rrvert \leqslant C \|g\|_p
\frac{\log|x|_*}{|x|_*^d}.
\end{equation}
This is true since $|\nabla_i \nabla_j G_{\mathsf{h}}(x)| \lesssim
|x|_*^{-d}$,
\[
\bigl\llvert \bigl\langle g \tarpas \phi_j(x+\mathbf{e}_i)
\bigr\rangle \bigr\rrvert \leqslant\|g\|_2 \tarpas \|\phi
_j\|_2,
\]
$\|\phi_j\|_2$ is finite by Theorem~\ref{tcorrector}, and we assume $p
\geqslant2$.

We now turn to the proof of (\ref{ez-two}). From Proposition~\ref
{prep-z}, we learn that
\[
\nabla z(e) = \sum_{y \in\mathbb{Z}^d} \nabla G(e,y) \bigl(R :
\nabla^2 G_{\mathsf{h}}\bigr) (y) + \sum
_{b
\in\mathbb{B}} \nabla\nabla G(e,b) h(b).
\]
We now proceed to show that each of the two terms
%
\begin{eqnarray}
\label{part1}
&& \sum_{y \in\mathbb{Z}^d} \bigl\llvert \bigl\langle
g \tarpas \nabla G(e,y) \bigl(R : \nabla^2 G_{\mathsf{h}}\bigr) (y)
\bigr\rangle \bigr\rrvert,
\\
\label{part2}
&& \sum_{b \in\mathbb{B}} \bigl\llvert \bigl\langle
g \tarpas \nabla\nabla G(e,b) h(b) \bigr\rangle \bigr\rrvert
\end{eqnarray}
is bounded by the right-hand side of (\ref{ez-two}).

\begin{longlist}[\textit{Step} I.1.]
\item[\textit{Step} I.1.]
We begin with (\ref{part1}), which is the
more delicate. As noted in Remark~\ref{rcentredR}, the random variable
$R$ is centered, so the expectation appearing within the absolute value
in (\ref{part1}) is in fact a correlation. We thus wish to apply
Corollary~\ref{ccovar} and write
%
\begin{eqnarray}
&& \bigl\llvert \bigl\langle g \nabla G(e,y) \bigl(R :
\nabla^2 G_{\mathsf
{h}}\bigr) (y) \bigr\rangle \bigr\rrvert
\nonumber
\\[-8pt]
\label{part1step1}
\\[-8pt]
\nonumber
&&\qquad\leqslant\sum_{b \in\mathbb{B}} \bigl\langle\bigl[
\partial_b \bigl(g \tarpas \nabla G(e,y)\bigr)\bigr]^2
\bigr\rangle^{1/2} \bigl\langle\bigl[\partial _b \bigl(R :
\nabla^2 G_{\mathsf{h}}\bigr) (y)\bigr]^2 \bigr\rangle^{1/2}.
\end{eqnarray}
However, recalling that
\[
(R - A_{\mathsf{h}})_{ij}(y) = - \bigl[A_i\bigl(
\mathbf{1}_j^i + \nabla _i
\phi_j\bigr) \bigr](y - \mathbf{e}_i),
\]
we see that a slight difficulty appears because we have not given a
meaning to $\partial_b \phi_j$. As was anticipated in Remark~\ref
{rphixmphimu}, this need not bother us. If we formally extend
Proposition~\ref{pderiv-corr} to the case $\mu= 0$, we arrive at the
formal expression
%
\begin{eqnarray}
\partial_b \bigl(R_{ij}(y)\bigr) &=&
\partial_b a_b \bigl(- \mathbf{1}_{b = (y-\mathbf{e}_i,y)}\bigl(
\mathbf {1}_j^i + \nabla_i \phi
_j\bigr) (y-\mathbf{e}_i)
\nonumber
\\[-8pt]
\label{formal1}
\\[-8pt]
\nonumber
&&{}+ A_i(y-\mathbf{e}_i) \nabla\nabla G(y-
\mathbf{e}_i,b) (\xi+ \nabla \phi _j) (b) \bigr).
\end{eqnarray}
The point now is that although we do not wish to discuss the sense of
(\ref{formal1}) as a derivative, we can take it as a \emph{definition}
of the random variable $\partial_b (R_{ij}(y))$, and observe that
(\ref{part1step1}) holds. To see this, we approximate the left-hand
side of (\ref{part1step1}) by introducing a small mass $\mu> 0$. We introduce
\[
\bigl(A_{\mathsf{h}}^\mu\bigr)_{ij} = \bigl\langle
A_i\bigl(\mathbf{1}_j^i +
\nabla_i \phi_{j,\mu}\bigr) \bigr\rangle
\]
and $R^{\mu}$ by setting
\[
\bigl(R^{\mu} - A_{\mathsf{h}}^{\mu}\bigr)_{ij}(y)
= - \bigl[A_i\bigl(\mathbf {1}_j^i +
\nabla_i \phi _{j,\mu
}\bigr) \bigr](y - \mathbf{e}_i)
\]
(where of course $\phi_{j,\mu} = \phi_{\mathbf{e}_j,\mu}$). We can
now write
the left-hand side of (\ref{part1step1}) as the limit as $\mu$ tends
to $0$ of
\[
\bigl\llvert \bigl\langle g \tarpas \nabla G(e,y) \bigl(R^\mu:
\nabla^2 G_{\mathsf
{h}}\bigr) (y) \bigr\rangle \bigr\rrvert.
\]
Applying Proposition~\ref{pHJ} on this term is now legitimate, and by
Proposition~\ref{pderiv-corr},
%
\begin{eqnarray}
 \partial_b \bigl(R^\mu_{ij}(y)
\bigr) &=& \partial_b a_b \bigl( - \mathbf{1}_{b = (y-\mathbf{e}_i,y)}
\bigl(\mathbf {1}_j^i + \nabla_i \phi
_{j,\mu
}\bigr) (y-\mathbf{e}_i)
\nonumber
\\[-8pt]
\label{non-formal1}
\\[-8pt]
\nonumber
&&{}+ A_i(y-\mathbf{e}_i) \nabla\nabla G_\mu(y-
\mathbf{e}_i,b) (\xi+ \nabla \phi_{j,\mu}) (b) \bigr).
\end{eqnarray}
By taking the limit $\mu\to0$, it follows that (\ref{part1step1})
holds with $\partial_b R$ defined by (\ref{formal1}).

\item[\textit{Step} I.2.]
By H\"{o}lder's inequality, it follows from Theorems~\ref{tcorrector}
and \ref{tmarott} that
\[
\bigl\langle\bigl(\partial_b R : \nabla^2G_{\mathsf{h}}
\bigr)^2(y) \bigr\rangle^{1/2} \lesssim\frac{1}{|\underline
{b}-y|_*^d \tarpas  |y|_*^d},
\]
where $\lesssim$ stands for $\leqslant$ up to a multiplicative
constant that
only depends on $d$ and the Lipschitz constant of $a$.
On the other hand, using Proposition~\ref{pderiv-green}, we see that
\begin{eqnarray*}
\partial_b \bigl(g \tarpas \nabla G(e,y)\bigr) & = & (
\partial_b g) \nabla G(e,y) + g \tarpas\, \partial_b
\nabla G(e,y)
\\
& = & (\partial_b g) \nabla G(e,y) - g \tarpas (
\partial_b a_b) \nabla \nabla G(e,b) \nabla G(y,b).
\end{eqnarray*}
Using H\"{o}lder's inequality (in conjunction with the strict
inequality $p>2$) and Theorem~\ref{tmarott}, we are led to
\[
\bigl\langle\bigl[\partial_b \bigl(g \tarpas \nabla G(e,y)\bigr)
\bigr]^2 \bigr\rangle^{1/2} \lesssim\frac{\| \partial_b g
\|
_p}{|\underline{e}-y|_*^{d-1}} +
\frac{\|g\|_p}{|\underline
{b}-\underline{e}|_*^d \tarpas  |\underline
{b}-y|_*^{d-1}}.
\]
So we obtain from (\ref{part1step1}) the inequality
%
\begin{eqnarray}
&& \bigl\llvert \bigl\langle g \tarpas \nabla G(e,y) \bigl(R :
\nabla^2 G_{\mathsf
{h}}\bigr) (y) \bigr\rangle \bigr\rrvert
\nonumber
\\[-8pt]
\label{eestim1}
\\[-8pt]
\nonumber
&&\qquad\lesssim\sum_{b \in\mathbb{B}} \biggl( \frac{\| \partial_b g \|
_p}{|\underline
{e}-y|_*^{d-1}} +
\frac{\|g\|_p}{|\underline{b}-\underline{e}|_*^d \tarpas  |\underline
{b}-y|_*^{d-1}} \biggr) \frac{1}{|\underline{b}-y|_*^d \tarpas  |y|_*^d},
\end{eqnarray}
and the term appearing in (\ref{part1}) is bounded (up to a constant) by
\[
\mathop{\sum_{y \in\mathbb{Z}^d}}_{b \in\mathbb{B}} \biggl(
\frac{\| \partial_b g \|
_p}{|\underline{e}-y|_*^{d-1}} + \frac{\|g\|_p}{|\underline
{b}-\underline{e}|_*^d \tarpas  |\underline
{b}-y|_*^{d-1}} \biggr) \frac{1}{|\underline{b}-y|_*^d \tarpas  |y|_*^d}.
\]
To see that this is bounded by the right-hand side of (\ref
{ez-two}), it
suffices to observe that
\[
\sum_{y \in\mathbb{Z}^d} \frac{1}{|\underline{b}-y|_*^{2d-1} \tarpas  |y|_*^d} \lesssim
\frac{1}{|\underline{b}|_*^d}
\]
and
%
\begin{equation}
\label{logcorr}
\sum_{b \in\mathbb{B}} \frac{1}{|\underline{b}-\underline{e}|_*^d
\tarpas  |\underline{b}|_*^d}
\lesssim \frac{\log|\underline{e}|_*}{|\underline{e}|_*^d}.
\end{equation}
These two facts are proved in Proposition~\ref{psums} of the \hyperref[sappend]{Appendix}.

\item[\textit{Step} II.]
We now turn to the analysis of (\ref{part2}).
We note that
\[
\sum_{b \in\mathbb{B}} \bigl\llvert \bigl\langle g \tarpas
\nabla\nabla G(e,b) h(b) \bigr\rangle \bigr\rrvert \leqslant\sum
_{b
\in\mathbb{B}} \|g\|_2 \bigl\langle \bigl( \nabla\nabla
G(e,b) h(b) \bigr)^2 \bigr\rangle^{1/2}.
\]
Using the explicit form of $h$ given by Proposition~\ref{pz-eq}
together with Theorems~\ref{tcorrector} and~\ref{tmarott}, we
arrive at
\[
\bigl\langle \bigl( \nabla\nabla G(e,b) h(b) \bigr)^2 \bigr
\rangle^{1/2} \lesssim\frac
{1}{|\underline
{b}-\underline{e}|_*^d \tarpas  |\underline{b}|_*^d}.
\]
In view of (\ref{logcorr}), we have shown that the term in (\ref
{part2}) is bounded by a constant times
\[
\|g\|_2 \tarpas \frac{\log|\underline{e}|_*}{|\underline{e}|_*^d},
\]
which is a better bound than needed.\quad\qed
\end{longlist}
\noqed\end{pf*}

\section{Proof of Theorem~\texorpdfstring{\protect\ref{tstruct}}{2.1}}
\label{sproof}
%
Our starting point is the identity
%
\begin{equation}
\label{start} \bigl\langle\phi_\xi(0) \phi_\xi(x) \bigr
\rangle = \sum_{e
\in\mathbb{B}} \bigl\langle\partial_e
\phi_\xi (0) \tarpas (\mathscr {L}+ 1)^{-1} \,\tarpas
\partial_e \phi_\xi(x) \bigr\rangle,
\end{equation}
with
%
\begin{equation}
\label{defdrphi} \partial_e \phi_\xi(y) = -
\partial_e a_e \nabla G(y,e) (\xi+ \nabla
\phi_\xi) (e) \qquad\bigl(y \in\mathbb{Z}^d\bigr).
\end{equation}
As in step~I.1 of the proof of Theorem~\ref{ttwo-scale}, we do not
mean to discuss the meaning of $\partial_e \phi_\xi(y)$ as a
derivative of
$\phi_\xi(y)$. Rather, it suffices for our purpose to observe that the
identity in (\ref{start}) holds with $\partial_e \phi_\xi(0)$ and
$\partial
_e \phi
_\xi(x)$ \emph{defined} by (\ref{defdrphi}). This follows easily by
approximating $\phi_\xi$ by $\phi_{\xi,\mu}$, applying
Propositions~\ref{pHJ}
and~\ref{pderiv-corr}, and letting $\mu$ tend to $0$.

Replacing $\partial_e \phi_\xi(0)$ and $\partial_e \phi_\xi(x)$
by their
definitions, the summand in the right-hand side of (\ref{start}) becomes
%
\begin{equation}
\label{esummand}
\hspace*{5pt}\bigl\langle\partial_e a_e \nabla
G(0,e) (\xi+ \nabla\phi_\xi ) (e) \tarpas (\mathscr{L}+
1)^{-1} \tarpas\, \partial_e a_e \nabla G(x,e)
(\xi+ \nabla\phi_\xi) (e) \bigr\rangle.
\end{equation}
We see that two $\nabla G$ terms appear in this expectation. We will
``pull out of the expectation'' each of these $\nabla G$ terms using
Theorem~\ref{ttwo-scale}. These form the two first steps of the proof.
The last step discusses how to replace $\nabla G_{\mathsf{h}}$ by its
continuous-space counterpart $\nabla\mathcal{G}_{\mathsf{h}}$.

\begin{longlist}[\textit{Step} 2.]
\item[\textit{Step} 1.] Defining
\[
g_e(x) = \partial_e a_e (\xi+ \nabla
\phi_\xi) (e) \tarpas (\mathscr{L}+ 1)^{-1} \tarpas
\,\partial_e a_e \nabla G(x,e) (\xi+ \nabla
\phi_\xi) (e),
\]
we see that we can rewrite the term in (\ref{esummand}) as
\[
\bigl\langle g_e(x) \tarpas \nabla G(0,e) \bigr\rangle,
\]
and we wish to justify that
%
\begin{equation}
\label{estep1todo}
\qquad\sum_{e \in\mathbb{B}} \Biggl\llvert \bigl\langle
g_e(x) \tarpas \nabla G(0,e) \bigr\rangle - \sum
_{j = 1}^d \nabla_j G_{\mathsf{h}}(e)
\bigl\langle g_e(x) \tarpas (\mathbf{e}_j + \nabla\phi
_j) (e) \bigr\rangle \Biggr\rrvert \lesssim\frac{\log
^2|x|_*}{|x|_*^{d-1}}.
\end{equation}
In order to apply Theorem~\ref{ttwo-scale} for this purpose, we need
to compute $\partial_b g_e(x)$ for every $b \in\mathbb{B}$. From the
commutation
relation in~(\ref{commutator}), it follows that
\[
\partial_b \tarpas \mathscr{L}= (\mathscr{L}+1) \tarpas
\,\partial_b,
\]
and thus
\[
\partial_b \tarpas (\mathscr{L}+1)^{-1} = (
\mathscr{L}+2)^{-1} \tarpas \,\partial_b.
\]
From this observation, we get that
%
\begin{eqnarray}
\label{edecompgex}
\partial_b g_e(x) &=&
g^{(1)}_{b,e}(x) + g^{(2)}_{b,e}(x) +
g^{(3)}_{b,e}(x) + g^{(4)}_{b,e}(x)
\end{eqnarray}
with
\begin{eqnarray*}
g^{(1)}_{b,e}(x) &=& - \partial_e a_e
\,\partial_b a_b \nabla\nabla G(e,b) (\xi+ \nabla \phi
_\xi) (b) \tarpas (\mathscr{L}+ 1)^{-1}\\
&&{}{}\times \tarpas
\,\partial_e a_e \nabla G(x,e) (\xi + \nabla
\phi_\xi) (e),
\\
g^{(2)}_{b,e}(x) &=& - \partial_e
a_e (\xi+ \nabla\phi_\xi) (e) \tarpas (\mathscr{L}+
2)^{-1}\partial_e a_e\,\partial_b a_b
\nabla\nabla G(e,b)
\\
&&{}\times\bigl[ \nabla G(x,b) (\xi+ \nabla\phi_\xi) (e)+
\nabla G(x,e) (\xi+ \nabla\phi_\xi) (b) \bigr],
\\
g^{(3)}_{b,e}(x) &=&  \mathbf{1}_{e = b} \tarpas
\,\partial_e^2 a_e (\xi+ \nabla
\phi_\xi) (e) \tarpas (\mathscr{L}+ 1)^{-1} \tarpas
\,\partial_e a_e \nabla G(x,e) (\xi+ \nabla
\phi_\xi) (e)
\end{eqnarray*}
and
\[
g^{(4)}_{b,e}(x) = \mathbf{1}_{e = b} \tarpas
\,\partial_e a_e (\xi+ \nabla \phi_\xi) (e)
\tarpas (\mathscr{L} +2)^{-1} \tarpas\, \partial_e^2
a_e \nabla G(x,e) (\xi+ \nabla\phi_\xi) (e).
\]
As before, we do not wish to discuss the meaning of (\ref{edecompgex})
as a derivative, but rather use the fact that if $\partial_b g_e(x)$ is
defined in this way, then by the usual approximation argument,
\begin{eqnarray*}
&& \Biggl\llvert \bigl\langle g_e(x) \tarpas \nabla G(0,e) \bigr
\rangle - \sum_{j
= 1}^d
\nabla_j G_{\mathsf{h}}(e) \bigl\langle g_e(x) \tarpas
(\mathbf{e}_j + \nabla\phi_j) (e) \bigr\rangle \Biggr
\rrvert
\\
&&\qquad\lesssim\bigl\|g_e(x)\bigr\|_p \frac{\log|\underline{e}|_*}{|\underline
{e}|_*^d} +
\mathop{\sum_{y \in\mathbb{Z}^d}}_{b \in\mathbb{B}} \bigl\llVert
\partial_b g_e(x) \bigr\rrVert _p
\frac
{1}{|\underline{e}-y|_*^{d-1} \tarpas  |\underline{b}-y|_*^{d} \tarpas  |y|_*^d}.
\end{eqnarray*}
From Proposition~\ref{pcontract} and Theorems~\ref{tcorrector} and
\ref{tmarott}, we learn that
\[
\bigl\|g_e(x)\bigr\|_p \lesssim\frac{1}{|\underline{e}-x|_*^{d-1}}
\]
and
\[
\bigl\|\partial_b g_e(x)\bigr\|_p \lesssim
\frac{1}{|\underline{b} -
\underline{e}|_*^d} \biggl(\frac
{1}{|\underline{e}-x|_*^{d-1}} + \frac{1}{|\underline
{b}-x|_*^{d-1}} \biggr).
\]
Hence, up to a multiplicative constant, the left-hand side of (\ref
{estep1todo}) is smaller than the sum of the following two terms:
%
\begin{eqnarray}
\label{eterm1}
&& \sum_{e \in\mathbb{B}} \frac{\log|\underline{e}|_*}{|\underline
{e}-x|_*^{d-1} \tarpas  |\underline
{e}|_*^{d}},
\\
\label{eterm2}
&& \mathop{\sum_{e,b \in\mathbb{B}}}_{y \in\mathbb{Z}^d}
\frac
{1}{|\underline{b} - \underline
{e}|_*^d} \biggl(\frac{1}{|\underline{e}-x|_*^{d-1}} + \frac
{1}{|\underline
{b}-x|_*^{d-1}} \biggr)
\frac{1}{|\underline{e}-y|_*^{d-1} \tarpas  |\underline{b}-y|_*^{d}
\tarpas  |y|_*^d}.
\end{eqnarray}
By Remark~\ref{rlogs} of the \hyperref[sappend]{Appendix}, the sum in (\ref{eterm1}) is
dominated by a constant times the right-hand side of (\ref
{estep1todo}). As for
the sum in (\ref{eterm2}), we can further split it into the sum of
%
\begin{equation}
\label{eterm21}
\mathop{\sum_{e,b \in\mathbb{B}}}_{y \in\mathbb{Z}^d}
\frac
{1}{|\underline{b} - \underline
{e}|_*^d \tarpas  |\underline{e}-x|_*^{d-1} \tarpas  |\underline{e}-y|_*^{d-1} \tarpas  |\underline
{b}-y|_*^{d} \tarpas  |y|_*^d}
\end{equation}
and
%
\begin{equation}
\label{eterm22}
\mathop{\sum_{e,b \in\mathbb{B}}}_{y \in\mathbb{Z}^d}
\frac
{1}{|\underline{b} - \underline
{e}|_*^d \tarpas  |\underline{b}-x|_*^{d-1} \tarpas  |\underline{e}-y|_*^{d-1} \tarpas  |\underline
{b}-y|_*^{d} \tarpas  |y|_*^d}.
\end{equation}
By repeatedly applying Proposition~\ref{psums} of the \hyperref[sappend]{Appendix}, we can
bound the sum in (\ref{eterm21}) by
\[
\mathop{\sum_{e\in\mathbb{B}}}_{y \in\mathbb{Z}^d}
\frac{\log
|\underline
{e}-y|_*}{|\underline
{e}-y|_*^{2d-1} \tarpas  |\underline{e}-x|_*^{d-1} \tarpas  |y|_*^d} \lesssim\sum_{e \in
\mathbb{B}
}
\frac{1}{|\underline{e}-x|_*^{d-1} \tarpas  |\underline{e}|_*^d} \lesssim \frac{\log
|x|_*}{|x|_*^{d-1}},
\]
and similarly, bound the sum in (\ref{eterm22}) by
\[
\mathop{\sum_{b \in\mathbb{B}}}_{y \in\mathbb{Z}^d}\frac{\log
|\underline
{b}-y|_*}{|\underline
{b}-y|_*^{2d-1} \tarpas  |\underline{b}-x|_*^{d-1} \tarpas  |y|_*^d}
\lesssim\sum_{b \in
\mathbb{B}
} \frac{1}{|\underline{b}|_*^d \tarpas  |\underline{b}-x|_*^{d-1} } \lesssim
\frac{\log
|x|_*}{|x|_*^{d-1}},
\]
and the proof of (\ref{estep1todo}) is complete.

\item[\textit{Step} 2.]
Recall that we have written $ \langle
\phi _\xi (0) \tarpas  \phi_\xi(x)  \rangle$ as
\[
\sum_{e \in\mathbb{B}} \bigl\langle g_e(x) \nabla
G(0,e) \bigr\rangle,
\]
so we proved in step~1 that
%
\begin{equation}
\label{estep1result}
\quad\Biggl\llvert \bigl\langle\phi_\xi(0) \tarpas
\phi_\xi(x) \bigr\rangle - \sum_{e \in\mathbb{B}} \sum
_{j = 1}^d \nabla_j
G_{\mathsf{h}}(e) \bigl\langle g_e(x) \tarpas (
\mathbf{e}_j + \nabla\phi_j) (e) \bigr\rangle \Biggr
\rrvert \lesssim\frac
{\log^2|x|_*}{|x|_*^{d-1}}.
\end{equation}
We now aim to show that
%
\begin{eqnarray}
&&\qquad \sum_{e \in\mathbb{B}} \sum
_{j = 1}^d \Biggl\llvert \nabla_j
G_{\mathsf
{h}}(e) \bigl\langle g_e(x) \tarpas (
\mathbf{e}_j + \nabla\phi_j) (e) \bigr\rangle - \sum
_{k = 1}^d \nabla_j
G_{\mathsf{h}}(e) \mathsf {Q}^{(\xi,e)}_{jk}
\nabla_k G_{\mathsf{h}}(e-x) \Biggr\rrvert
\nonumber
\\[-8pt]
\label{eaim2}
\\[-8pt]
\nonumber
&&\qquad\qquad\lesssim\frac{\log^2|x|_*}{|x|_*^{d-1}},
\end{eqnarray}
where $\mathsf{Q}^{(\xi,e)}_{jk}$ is defined by
%
\begin{eqnarray}
\mathsf{Q}^{(\xi,e)}_{jk}&=& \bigl\langle\partial_e a_e (
\mathbf{e}_j + \nabla\phi_j) (e) (\xi + \nabla\phi
_\xi) (e)
\nonumber
\\[-8pt]
\label{defQxe}
\\[-8pt]
\nonumber
&&{}\times (\mathscr{L}+ 1)^{-1} \tarpas
\,\partial_e a_e (\mathbf{e}_k + \nabla
\phi_k) (e) (\xi+ \nabla\phi _\xi) (e) \bigr\rangle.
\end{eqnarray}
For $j \in\{1,\ldots,d\}$, we let
\[
\tilde{g}_{e,j} = \partial_e a_e (\xi+
\nabla\phi_\xi) (e) \tarpas (\mathscr{L}+ 1)^{-1} \tarpas
\,\partial _e a_e (\mathbf{e}_j + \nabla
\phi_j) (e) (\xi+ \nabla\phi_\xi) (e),
\]
and observe that since $(\mathscr{L}+1)^{-1}$ is symmetric,
%
\begin{equation}
\label{eselfadj}
\bigl\langle g_e(x) \tarpas (\mathbf{e}_j
+ \nabla\phi_j) (e) \bigr\rangle = \bigl\langle\tilde
{g}_{e,j} \tarpas \nabla G(x,e) \bigr\rangle.
\end{equation}
We let
\[
\partial_b \tilde{g}_{e,j} = \tilde{g}^{(1)}_{b,e,j}
+ \tilde {g}^{(1)}_{b,e,j} + \tilde {g}^{(3)}_{b,e,j}
+ \tilde{g}^{(4)}_{b,e,j},
\]
where
\begin{eqnarray*}
\tilde{g}^{(1)}_{b,e,j} &=& - \partial_e
a_e\, \partial_b a_b \nabla \nabla G(e,b) (
\xi+ \nabla \phi_\xi) (b) \tarpas (\mathscr{L}+ 1)^{-1}
\tarpas\\
&&{}\times \partial_e a_e (\mathbf {e}_j +
\nabla\phi _j) (e) (\xi+ \nabla\phi_\xi) (e),
\\
\tilde{g}^{(2)}_{b,e,j} &=& -\partial_e
a_e (\xi+ \nabla\phi_\xi) (e) \tarpas (\mathscr{L}+
2)^{-1}
\\
&&{}\times\partial_e a_e\,\partial_b a_b
\nabla\nabla G(e,b) \\
&&{}\times\bigl[(\mathbf{e}_j + \nabla \phi_j)
(b) (\xi+ \nabla\phi_\xi) (e) + (\mathbf{e}_j + \nabla
\phi_j) (e) (\xi + \nabla \phi_\xi) (b) \bigr],
\\
\tilde{g}^{(3)}_{b,e,j} &=& \mathbf{1}_{e=b} \tarpas
\,\partial_e^2 a_e (\xi+ \nabla
\phi_\xi ) (e) \tarpas (\mathscr{L}+ 1)^{-1} \tarpas
\,\partial_e a_e (\mathbf{e}_j + \nabla\phi
_j) (e) (\xi+ \nabla \phi_\xi) (e)
\end{eqnarray*}
and
\[
\tilde{g}^{(4)}_{b,e,j} = \mathbf{1}_{e=b} \tarpas
\,\partial_e a_e (\xi+ \nabla\phi_\xi) (e)
\tarpas (\mathscr{L}+ 2)^{-1} \tarpas \,\partial_e^2
a_e (\mathbf{e}_j + \nabla\phi _j) (e) (
\xi+ \nabla\phi_\xi) (e).
\]
As before (and because of Remark~\ref{rtranslation}), this definition
ensures that
\begin{eqnarray*}
&&\Biggl\llvert \bigl\langle\tilde{g}_{e,j} \tarpas \nabla G(x,e) \bigr
\rangle - \sum_{k = 1}^d
\nabla_k G_{\mathsf{h}}(e-x) \bigl\langle\tilde{g}_{e,j}
\tarpas (\mathbf{e}_k + \nabla\phi_k) (e) \bigr\rangle
\Biggr\rrvert
\\
&&\qquad\lesssim\|\tilde{g}_{e,j}\|_p \frac{\log|\underline
{e}-x|_*}{|\underline
{e}-x|_*^d} +
\mathop{\sum_{y \in\mathbb{Z}^d}}_{b \in\mathbb{B}} \llVert
\partial_b \tilde {g}_{e,j} \rrVert _p
\frac{1}{|\underline{e}-y|_*^{d-1} \tarpas  |\underline{b}-y|_*^{d} \tarpas  |y-x|_*^d}.
\end{eqnarray*}
Moreover, we infer from Proposition~\ref{pcontract} and Theorems~\ref
{tcorrector} and \ref{tmarott} that for any $1\leqslant p<\infty$
(and thus in particular the $p>2$ needed above)
\[
\|\tilde{g}_{e,j}\|_p \lesssim1
\]
and
\[
\llVert \partial_b \tilde{g}_{e,j} \rrVert _p
\lesssim\frac
{1}{|\underline{b}-\underline{e}|_*^d}.
\]
Since
\[
\bigl\langle\tilde{g}_{e,j} \tarpas (\mathbf{e}_k +
\nabla\phi_k) (e) \bigr\rangle = \mathsf{Q}^{(\xi,e)}_{jk},
\]
we obtain that
\begin{eqnarray*}
&& \Biggl\llvert \bigl\langle\tilde{g}_{e,j} \tarpas \nabla
G(x,e) \bigr\rangle - \sum_{k = 1}^d
\mathsf{Q}^{(\xi,e)}_{jk} \nabla _k
G_{\mathsf{h}}(e-x) \Biggr\rrvert
\\
&&\qquad\lesssim\frac{\log|\underline{e}-x|_*}{|\underline{e}-x|_*^d} + \mathop{\sum_{y \in\mathbb{Z}^d}}_{b \in\mathbb{B}}
\frac{1}{|\underline
{b}-\underline{e}|_*^d \tarpas  |\underline
{e}-y|_*^{d-1} \tarpas  |\underline{b}-y|_*^{d} \tarpas  |y-x|_*^d},
\end{eqnarray*}
and thus by (\ref{eselfadj}), up to a multiplicative constant, the
left-hand side
of (\ref{eaim2}) is smaller than
%
\begin{equation}
\label{erema}
\qquad\sum_{e \in\mathbb{B}} \frac{1}{|e|_*^{d-1}} \biggl(
\frac{\log
|\underline
{e}-x|_*}{|\underline
{e}-x|_*^d} + \mathop{\sum_{y \in\mathbb{Z}^d}}_{b \in\mathbb
{B}}
\frac{1}{|\underline{b}-\underline{e}|_*^d \tarpas  |\underline{e}-y|_*^{d-1} \tarpas  |\underline
{b}-y|_*^{d} \tarpas  |y-x|_*^d} \biggr).
\end{equation}
From Remark~\ref{rlogs} of the \hyperref[sappend]{Appendix}, we have
\[
\sum_{e \in\mathbb{B}} \frac{1}{|e|_*^{d-1}} \frac{\log
|\underline
{e}-x|_*}{|\underline
{e}-x|_*^d}
\lesssim\frac{\log^2|x|_*}{|x|_*^{d-1}}.
\]
The remaining sum from (\ref{erema}) can be bounded, using
Proposition~\ref{psums} repeatedly, by
\[
\mathop{\sum_{y \in\mathbb{Z}^d}}_{e \in\mathbb{B}}
\frac
{\log|\underline
{e}-y|_*}{|\underline
{e}|_*^{d-1} \tarpas   |\underline{e}-y|_*^{2d-1} \tarpas  |y-x|_*^d} \lesssim\sum_{y \in
\mathbb{Z}
^d}
\frac{1}{|y|_*^{d-1} \tarpas  |y-x|_*^d} \lesssim\frac{\log|x|_*}{|x|_*^{d-1}},
\]
and this finishes the proof of (\ref{eaim2}).

\item[\textit{Step} 3.] Note that by the stationarity of the
environment, the matrix $\mathsf{Q}^{(\xi,e)}$ depends on the edge
$e$ only through its
orientation. On the other hand, the quantities $\nabla_j G_{\mathsf
{h}}(e)$ and
$\nabla
_j G_{\mathsf{h}}(e-x)$ depend on the edge $e$ only through its base point.
We also observe that the matrix $\mathsf{Q}^{(\xi)}$ introduced in
(\ref{defQx}) is by
definition $\sum_{e \in\mathcal{E}_0} \mathsf{Q}^{(\xi,e)}$.
Hence, the
previous steps of the proof have led us  [see (\ref{estep1result}) and
(\ref{eaim2})] to
%
\begin{equation}
\label{estep2result} \hspace*{8pt}\Biggl\llvert \bigl\langle\phi_\xi(0) \tarpas
\phi_\xi(x) \bigr\rangle - \sum_{y \in\mathbb{Z}^d} \sum
_{j,k =
1}^d \nabla_j
G_{\mathsf{h}}(y) \mathsf{Q}^{(\xi)}_{jk} \nabla_k
G_{\mathsf{h}}(y-x) \Biggr\rrvert \lesssim\frac
{\log
^2|x|_*}{|x|_*^{d-1}}.
\end{equation}
In order to complete the proof of Theorem~\ref{tstruct}, it thus
suffices to show that
\[
\Biggl\llvert \sum_{y \in\mathbb{Z}^d} \sum
_{j,k = 1}^d \nabla_j G_{\mathsf{h}}(y)
\mathsf{Q}^{(\xi)}_{jk} \nabla_k G_{\mathsf{h}}
(y-x) - \mathscr{K}_\xi(x) \Biggr\rrvert \lesssim\frac{\log|x|_*}{|x|_*^{d-1}},
\]
where $\mathscr{K}_\xi$ was introduced in (\ref{defK}). We learn from
Proposition~\ref{pdiscrete-to-continuous} of the \hyperref[sappend]{Appendix} that
\[
\biggl\llvert \nabla_j G_{\mathsf{h}}(y) - \frac{\partial\mathcal
{G}_{\mathsf{h}}}{\partial y_j}(y)
\biggr\rrvert \lesssim \frac{1}{|y|^d}.
\]
As a consequence,
\begin{eqnarray*}
&&\mathop{\sum_{ y \in\mathbb{Z}^d \setminus\{0\}}}_{ 1 \leqslant
j,k \leqslant d} \biggl\llvert
\nabla_j G_{\mathsf{h}}(y) \mathsf{Q}^{(\xi)}_{jk}
\nabla_k G_{\mathsf{h}}(y-x) - \frac{\partial\mathcal{G}_{\mathsf
{h}}}{\partial
y_j}(y)
\mathsf{Q}^{(\xi)} _{jk} \nabla_k
G_{\mathsf{h}}(y-x) \biggr\rrvert
\\
&&\qquad\lesssim\sum_{y \in\mathbb{Z}^d \setminus\{0\}} \frac{1}{|y|^d \tarpas  |y-x|_*^{d-1}} \lesssim
\frac{\log|x|_*}{|x|_*^{d-1}},
\end{eqnarray*}
where we used Proposition~\ref{psums} of the \hyperref[sappend]{Appendix} in the last
step. Similarly,
\[
\mathop{\sum_{y \in\mathbb{Z}^d \setminus\{0,x\}}}_{1 \leqslant
j,k \leqslant d} \biggl\llvert
\frac{\partial\mathcal{G}_{\mathsf{h}}}{\partial y_j}(y) \mathsf {Q}^{(\xi)}_{jk}
\nabla_k G_{\mathsf{h}}(y-x) - \frac{\partial
\mathcal{G}_{\mathsf{h}}
}{\partial y_j}(y)
\mathsf{Q}^{(\xi)}_{jk} \frac{\partial\mathcal
{G}_{\mathsf{h}}}{\partial y_k}(y-x) \biggr\rrvert
\lesssim\frac{\log|x|_*}{|x|_*^{d-1}}.
\]
Moreover, one can check that
\begin{eqnarray*}
&&\mathop{\sum_{y \in\mathbb{Z}^d \setminus\{0,x\}}}_{1 \leqslant
j,k \leqslant d} \biggl\llvert
\frac{\partial\mathcal{G}_{\mathsf{h}}}{\partial y_j}(y)\mathsf {Q}^{(\xi)}_{jk}
\frac{\partial\mathcal{G}_{\mathsf{h}}}{\partial y_k}(y-x) - \int_{y+[0,1]^d} \frac{\partial\mathcal{G}_{\mathsf{h}}}{\partial
y_j}
\bigl(y'\bigr)\mathsf{Q}^{(\xi)}_{jk}
\frac{\partial\mathcal{G}_{\mathsf{h}}
}{\partial
y_k}\bigl(y'-x\bigr) \tarpas \,\mathrm{d}y'
\biggr\rrvert \\
&&\qquad\lesssim\frac{\log|x|_*}{|x|_*^{d-1}}.
\end{eqnarray*}
In these computations, we have been forced to drop some terms indexed
by $y \in\{0,x\}$. But it is easy to check that these terms are
negligible, for example,
\[
\mathop{\sum_{y \in\{0,x\}}}_{ 1 \leqslant j,k \leqslant d} \biggl\llvert
\int_{y+[0,1]^d} \frac{\partial\mathcal{G}_{\mathsf{h}}}{\partial
y_j}(u)\mathsf{Q}^{(\xi)}_{jk}
\frac{\partial\mathcal{G}_{\mathsf{h}}
}{\partial
y_k}(u-x) \tarpas \,{\mathrm{d}}u \biggr\rrvert \lesssim
\frac{1}{|x|_*^{d-1}} \qquad\bigl(x \in\mathbb{Z}^d \setminus\{0\}\bigr),
\]
so the proof is complete.
\end{longlist}

\begin{appendix}
\section*{Appendix: Basic estimates on discrete convolutions and Green functions}
\label{sappend}

\begin{prop}
\label{psums}
For every $\alpha> d$ and $\beta\in(0,\alpha]$,
\[
\sum_{y \in\mathbb{Z}^d} \frac{1}{|y|_*^\alpha\tarpas  |y-x|_*^\beta} \lesssim
\frac
{1}{|x|_*^\beta},
\]
while for $\beta\in(0,d]$,
\[
\sum_{y \in\mathbb{Z}^d} \frac{1}{|y|_*^d \tarpas  |y-x|_*^\beta} \lesssim
\frac
{\log
|x|_*}{|x|_*^\beta}.
\]
(In both statements, the sign $\lesssim$ hides a multiplicative
constant that does not depend on $x \in\mathbb{Z}^d$.)
\end{prop}
\begin{pf}
We give a unified proof of these two results, although it will be
apparent that the proof of the first statement alone can be slightly
simplified. We thus assume $\alpha\geqslant d$ and $\beta\in
(0,\alpha]$.
We decompose the sum over $y \in\mathbb{Z}^d$ according to whether
$|y| \geqslant
2|x|$ or not. If $|y| \geqslant2|x|$, then $|y-x| \geqslant|y|/2$,
and thus
\[
\sum_{|y| \geqslant2|x|} \frac{1}{|y|_*^\alpha\tarpas  |y-x|_*^\beta} \lesssim \sum
_{|y| \geqslant2|x|} \frac{1}{|y|_*^{\alpha+\beta}} \lesssim\frac
{1}{|x|_*^{\alpha+ \beta- d}}
\leqslant\frac{1}{|x|_*^{\beta}}
\]
(here and below, we understand that $y$ is the variable of summation).
We split the rest of the sum into two parts along the condition $|y-x|
\geqslant|x|/2$. This gives us two contributions, the first of which is
\[
\mathop{\sum_{|y| \leqslant2 |x|}}_{ |y-x| \geqslant|x|/2}
\frac
{1}{|y|_*^\alpha\tarpas  |y-x|_*^\beta} \lesssim\frac{1}{|x|_*^\beta} \sum_{|y| \leqslant2 |x|}
\frac{1}{|y|_*^\alpha}.
\]
This last sum is uniformly bounded if $\alpha> d$, while it is bounded
by $\log|x|_*$ if $\alpha= d$. For the second contribution to be
considered, note that $|y-x| \leqslant|x|/2$ implies that $|y|
\geqslant|x|/2$,
and thus
\[
\mathop{\sum_{|y| \leqslant2 |x|}}_{ |y-x| \leqslant|x|/2}
\frac
{1}{|y|_*^\alpha\tarpas  |y-x|_*^\beta} \lesssim\frac{1}{|x|_*^\alpha} \sum_{|y-x| \leqslant|x|/2}
\frac{1}{|y-x|_*^\beta}.
\]
Up to a constant, this last sum is bounded by
\[
\left|
\begin{array}{l@{\qquad}l}
1, & \mbox{if } \beta> d,
\cr\vspace*{2pt}
\log|x|_*,  & \mbox{if } \beta= d,
\cr\vspace*{4pt}
|x|_*^{d-\beta}, & \mbox{if } \beta< d.
\end{array}\right.
\]
Thus, this second contribution is always at most of the order of the
first, and this completes the proof.
\end{pf}
%
\begin{rem}
\label{rlogs}
The proof of Proposition~\ref{psums} can be adapted to yield, for
every $\beta\in(0,d]$,
\[
\sum_{y \in\mathbb{Z}^d} \frac{\log|y|_*}{|y|_*^d \tarpas  |y-x|_*^\beta} \lesssim
\frac{\log^2|x|_*}{|x|_*^\beta}.
\]
\end{rem}

\begin{prop}\label{pdiscrete-to-continuous}
For every $k \in\{1,\ldots, d \}$,
\[
\biggl\llvert \nabla_k G_{\mathsf{h}}(x) - \frac{\partial}{\partial x_k}
\mathcal{G}_{\mathsf{h}}(x) \biggr\rrvert \lesssim \frac{1}{|x|^d}.
\]
\end{prop}
\begin{pf}
Recall that $A_{\mathsf{h}}$ is a diagonal matrix, the diagonal
entries of which
we denote by $A_{\mathsf{h},1}, \ldots, A_{\mathsf{h},d}$. For $p
\in
[-\pi,\pi]^d$, let
\[
s(p) = 2 \sum_{j = 1}^d A_{\mathsf{h},j}
\bigl(1-\cos(p_j)\bigr).
\]
Using Fourier transforms, one can represent the Green function
$G_{\mathsf{h}}$ as
\[
G_{\mathsf{h}}(x) = \frac{1}{(2\pi)^d} \int_{\mathbb{T}}
\frac
{e^{-i p \cdot x}
}{s(p)} \tarpas \,{\mathrm{d}}p,
\]
where $\mathbb{T}= [-\pi,\pi]^d$. Similarly,
\[
\nabla_j G_{\mathsf{h}}(x) = \frac{1}{(2\pi)^d} \int
_{\mathbb{T}} \frac{(e^{-i p_j} - 1)
}{s(p)} \tarpas e^{-i p \cdot x} \tarpas \,{
\mathrm{d}}p.
\]
Let $\eta(x) = (2\pi)^{-d/2} e^{-|x|^2/2}$. We note that
\[
\biggl\llvert \frac{\partial\mathcal{G}_{\mathsf{h}}}{\partial x_k} (x) - \biggl(\frac{\partial\mathcal{G}_{\mathsf{h}}}{\partial x_k} * \eta \biggr)
(x) \biggr\rrvert \lesssim\frac{1}{|x|^d},
\]
where $*$ denotes the convolution. This can be seen, for instance, using
the explicit formula for the Green function,
\[
\mathcal{G}_{\mathsf{h}}(x) = \frac{1}{(d-2) \gamma_d \tarpas  |\det
(A_{\mathsf{h}})| \tarpas  (x \cdot A_{\mathsf{h}}^{-1}
x)^{(d-2)/2}},
\]
where $\gamma_d$ denotes the area measure of the unit sphere. The
regularization by convolution permits us to write down the Fourier
representation
\[
\biggl(\frac{\partial\mathcal{G}_{\mathsf{h}}}{\partial x_k} * \eta \biggr) (x) = \frac{1}{(2\pi)^d} \int
_{\mathbb{R}^d} - \frac{i p_j }{p \cdot A_{\mathsf{h}}p} \tarpas e^{-|p|^2/2}\tarpas
e^{-i p \cdot x} \tarpas \,{\mathrm{d}}p.
\]
In order to prove the proposition, it thus suffices to show that
\[
\biggl\llvert \int_{\mathbb{T}} \frac{(e^{-i p_j} - 1) }{s(p)} \tarpas
e^{-i p
\cdot x} \tarpas \,{\mathrm{d}}p - \int_{\mathbb{R}^d} -
\frac{i p_j }{p \cdot A_{\mathsf{h}}p} \tarpas e^{-|p|^2/2}\tarpas e^{-i p
\cdot x} \tarpas \,{
\mathrm{d}}p \biggr\rrvert \lesssim\frac{1}{|x|^d}.
\]
We select a smooth cut-off function $\chi(p)$ that is equal to one near
$p=0$ and is compactly supported in $\mathbb{T}$.
We use it to split the left-hand side into
\[
\int_{\mathbb{T}}(1-\chi) (p)\frac{(e^{-i p_j} - 1) }{s(p)} \tarpas
e^{-i
p \cdot
x} \,\tarpas {\mathrm{d}}p\quad\mbox{and}\quad \int
_{\mathbb{R}^d}f(p)\tarpas e^{-i p \cdot x} \tarpas \,{\mathrm{d}}p,
\]
where
\[
f(p) = \chi(p)\frac{(e^{-i p_j} - 1) }{s(p)} + \frac{i p_j }{p \cdot
A_{\mathsf{h}}p} \tarpas e^{-|p|^2/2}
\]
can be considered to be defined on all $\mathbb{R}^d$. By the
properties of $\chi$,
$(1-\chi)(p)\frac{(e^{-i p_j} - 1) }{s(p)}$ is a smooth periodic
function on $\mathbb{T}$, so that we obtain by
integrations by parts that
\[
\int_{\mathbb{T}} (1-\chi) (p)\frac{(e^{-i p_j} - 1) }{s(p)}\tarpas
e^{-i
p \cdot
x} \tarpas \,{\mathrm{d}}p
\]
decays faster than any negative power of $|x|$. Hence, it suffices to
show that
%
\setcounter{equation}{0}
\begin{equation}
\label{etoshowgreen}
\biggl\llvert \int_{\mathbb{R}^d} f(p) \tarpas
e^{-ip \cdot x} \tarpas \,{\mathrm{d}}p \biggr\rrvert \lesssim \frac{1}{|x|^d}.
\end{equation}
One can decompose $f$ as
\[
f(p) = -\frac{p_j^2}{2 p \cdot A_{\mathsf{h}}p} + \tilde{f}(p),
\]
so that $\tilde{f}$ is ``more regular'' than $f$ close to the origin. One
can then show by integration by parts that
\[
\biggl\llvert \int_{\mathbb{R}^d} -\frac{p_j^2}{2 p \cdot A_{\mathsf
{h}}p}\tarpas
e^{-ip
\cdot
x} \tarpas \,{\mathrm{d}}p - \biggl(\frac{\partial^2 \mathcal{G}_{\mathsf
{h}}}{\partial x_j^2}* \eta \biggr)
(x)\biggr\rrvert \lesssim \frac{1}{|x|^d}
\]
and
\[
\biggl\llvert \int_{\mathbb{R}^d} \tilde{f}(p) \tarpas
e^{-ip \cdot x} \tarpas \,{\mathrm{d}}p \biggr\rrvert \lesssim\frac{1}{|x|^d}
\]
for any $x\in\mathbb{R}^d$.
Since $(\partial^2 \mathcal{G}_{\mathsf{h}}/\partial x_j^2 * \eta
)(x) \lesssim|x|^{-d}$, the proof
is complete.
\end{pf}
\end{appendix}

\section*{Acknowledgement}
We would like to thank Marek
Biskup for stimulating discussions about this problem.






\printaddresses

\begin{thebibliography}{40}


\bibitem{armsma}
\begin{barticle}[mr]
\bauthor{\bsnm{Armstrong},~\bfnm{Scott~N.}\binits{S.~N.}} \AND
\bauthor{\bsnm{Smart},~\bfnm{Charles~K.}\binits{C.~K.}}
(\byear{2014}).
\btitle{Quantitative stochastic homogenization of elliptic equations in nondivergence form}.
\bjournal{Arch. Ration. Mech. Anal.}
\bvolume{214}
\bpages{867--911}.
\bid{doi={10.1007/s00205-014-0765-6}, issn={0003-9527}, mr={3269637}}
\end{barticle}
%
\bptok{imsref}%
\endbibitem

\bibitem{bal1}
\begin{barticle}[mr]
\bauthor{\bsnm{Bal},~\bfnm{Guillaume}\binits{G.}}
(\byear{2008}).
\btitle{Central limits and homogenization in random media}.
\bjournal{Multiscale Model. Simul.}
\bvolume{7}
\bpages{677--702}.
\bid{doi={10.1137/070709311}, issn={1540-3459}, mr={2443008}}
\end{barticle}
%
\bptok{imsref}%
\endbibitem

\bibitem{bal2}
\begin{barticle}[mr]
\bauthor{\bsnm{Bal},~\bfnm{Guillaume}\binits{G.}}
(\byear{2010}).
\btitle{Homogenization with large spatial random potential}.
\bjournal{Multiscale Model. Simul.}
\bvolume{8}
\bpages{1484--1510}.
\bid{doi={10.1137/090754066}, issn={1540-3459}, mr={2718269}}
\end{barticle}
%
\bptok{imsref}%
\endbibitem

\bibitem{bal3}
\begin{barticle}[mr]
\bauthor{\bsnm{Bal},~\bfnm{Guillaume}\binits{G.}}
(\byear{2011}).
\btitle{Convergence to homogenized or stochastic partial differential equations}.
\bjournal{Appl. Math. Res. Express. AMRX}
\bvolume{2}
\bpages{215--241}.
\bid{issn={1687-1200}, mr={2835990}}
\end{barticle}
%
\bptok{imsref}%
\endbibitem

\bibitem{balgar2}
\begin{barticle}[mr]
\bauthor{\bsnm{Bal},~\bfnm{Guillaume}\binits{G.}},
\bauthor{\bsnm{Garnier},~\bfnm{Josselin}\binits{J.}},
\bauthor{\bsnm{Gu},~\bfnm{Yu}\binits{Y.}} \AND
\bauthor{\bsnm{Jing},~\bfnm{Wenjia}\binits{W.}}
(\byear{2012}).
\btitle{Corrector theory for elliptic equations with long-range correlated random potential}.
\bjournal{Asymptot. Anal.}
\bvolume{77}
\bpages{123--145}.
\bid{issn={0921-7134}, mr={2977330}}
\end{barticle}
%
\bptok{imsref}%
\endbibitem

\bibitem{balgar}
\begin{barticle}[mr]
\bauthor{\bsnm{Bal},~\bfnm{Guillaume}\binits{G.}},
\bauthor{\bsnm{Garnier},~\bfnm{Josselin}\binits{J.}},
\bauthor{\bsnm{Motsch},~\bfnm{S{\'e}bastien}\binits{S.}} \AND
\bauthor{\bsnm{Perrier},~\bfnm{Vincent}\binits{V.}}
(\byear{2008}).
\btitle{Random integrals and correctors in homogenization}.
\bjournal{Asymptot. Anal.}
\bvolume{59}
\bpages{1--26}.
\bid{issn={0921-7134}, mr={2435670}}
\end{barticle}
%
\bptok{imsref}%
\endbibitem

\bibitem{balgu-rev}
\begin{barticle}[mr]
\bauthor{\bsnm{Bal},~\bfnm{Guillaume}\binits{G.}} \AND
\bauthor{\bsnm{Gu},~\bfnm{Yu}\binits{Y.}}
(\byear{2015}).
\btitle{Limiting models for equations with large random potential: A~review}.
\bjournal{Commun. Math. Sci.}
\bvolume{13}
\bpages{729--748}.
\bid{doi={10.4310/CMS.2015.v13.n3.a7}, issn={1539-6746}, mr={3318383}}
\bptnote{check pages, check year}%
\end{barticle}
%
\bptok{imsref}%
\endbibitem

\bibitem{baljin}
\begin{barticle}[mr]
\bauthor{\bsnm{Bal},~\bfnm{Guillaume}\binits{G.}} \AND
\bauthor{\bsnm{Jing},~\bfnm{Wenjia}\binits{W.}}
(\byear{2011}).
\btitle{Corrector theory for elliptic equations in random media with singular {G}reen's function. {A}pplication to random boundaries}.
\bjournal{Commun. Math. Sci.}
\bvolume{9}
\bpages{383--411}.
\bid{issn={1539-6746}, mr={2815677}}
\end{barticle}
%
\bptok{imsref}%
\endbibitem

\bibitem{berbis}
\begin{barticle}[mr]
\bauthor{\bsnm{Berger},~\bfnm{Noam}\binits{N.}} \AND
\bauthor{\bsnm{Biskup},~\bfnm{Marek}\binits{M.}}
(\byear{2007}).
\btitle{Quenched invariance principle for simple random walk on percolation clusters}.
\bjournal{Probab. Theory Related Fields}
\bvolume{137}
\bpages{83--120}.
\bid{doi={10.1007/s00440-006-0498-z}, issn={0178-8051}, mr={2278453}}
\end{barticle}
%
\bptok{imsref}%
\endbibitem

\bibitem{bisspo}
\begin{barticle}[mr]
\bauthor{\bsnm{Biskup},~\bfnm{Marek}\binits{M.}} \AND
\bauthor{\bsnm{Spohn},~\bfnm{Herbert}\binits{H.}}
(\byear{2011}).
\btitle{Scaling limit for a class of gradient fields with nonconvex potentials}.
\bjournal{Ann. Probab.}
\bvolume{39}
\bpages{224--251}.
\bid{doi={10.1214/10-AOP548}, issn={0091-1798}, mr={2778801}}
\end{barticle}
%
\bptok{imsref}%
\endbibitem

\bibitem{boivin}
\begin{barticle}[mr]
\bauthor{\bsnm{Boivin},~\bfnm{Daniel}\binits{D.}}
(\byear{2009}).
\btitle{Tail estimates for homogenization theorems in random media}.
\bjournal{ESAIM Probab. Stat.}
\bvolume{13}
\bpages{51--69}.
\bid{doi={10.1051/ps:2007036}, issn={1292-8100}, mr={2493855}}
\end{barticle}
%
\bptok{imsref}%
\endbibitem

\bibitem{boupia}
\begin{barticle}[mr]
\bauthor{\bsnm{Bourgeat},~\bfnm{Alain}\binits{A.}} \AND
\bauthor{\bsnm{Piatnitski},~\bfnm{Andrey}\binits{A.}}
(\byear{2004}).
\btitle{Approximations of effective coefficients in stochastic homogenization}.
\bjournal{Ann. Inst. Henri Poincar\'e Probab. Stat.}
\bvolume{40}
\bpages{153--165}.
\bid{doi={10.1016/S0246-0203(03)00065-7}, issn={0246-0203}, mr={2044813}}
\end{barticle}
%
\bptok{imsref}%
\endbibitem

\bibitem{cafsou}
\begin{barticle}[mr]
\bauthor{\bsnm{Caffarelli},~\bfnm{Luis~A.}\binits{L.~A.}} \AND
\bauthor{\bsnm{Souganidis},~\bfnm{Panagiotis~E.}\binits{P.~E.}}
(\byear{2010}).
\btitle{Rates of convergence for the homogenization of fully nonlinear uniformly elliptic pde in random media}.
\bjournal{Invent. Math.}
\bvolume{180}
\bpages{301--360}.
\bid{doi={10.1007/s00222-009-0230-6}, issn={0020-9910}, mr={2609244}}
\end{barticle}
%
\bptok{imsref}%
\endbibitem

\bibitem{capiof}
\begin{barticle}[mr]
\bauthor{\bsnm{Caputo},~\bfnm{Pietro}\binits{P.}} \AND
\bauthor{\bsnm{Ioffe},~\bfnm{Dmitry}\binits{D.}}
(\byear{2003}).
\btitle{Finite volume approximation of the effective diffusion matrix: The case of independent bond disorder}.
\bjournal{Ann. Inst. Henri Poincar\'e Probab. Stat.}
\bvolume{39}
\bpages{505--525}.
\bid{doi={10.1016/S0246-0203(02)00016-X}, issn={0246-0203}, mr={1978989}}
\end{barticle}
%
\bptok{imsref}%
\endbibitem

\bibitem{confah}
\begin{barticle}[mr]
\bauthor{\bsnm{Conlon},~\bfnm{Joseph~G.}\binits{J.~G.}} \AND
\bauthor{\bsnm{Fahim},~\bfnm{Arash}\binits{A.}}
(\byear{2015}).
\btitle{Strong convergence to the homogenized limit of parabolic equations with random coefficients}.
\bjournal{Trans. Amer. Math. Soc.}
\bvolume{367}
\bpages{3041--3093}.
\bid{doi={10.1090/S0002-9947-2014-06005-4}, issn={0002-9947}, mr={3314801}}
\bptnote{check volume, check pages}%
\end{barticle}
%
\bptok{imsref}%
\endbibitem

\bibitem{conspe}
\begin{barticle}[mr]
\bauthor{\bsnm{Conlon},~\bfnm{Joseph~G.}\binits{J.~G.}} \AND
\bauthor{\bsnm{Spencer},~\bfnm{Thomas}\binits{T.}}
(\byear{2014}).
\btitle{Strong convergence to the homogenized limit of elliptic equations with random coefficients}.
\bjournal{Trans. Amer. Math. Soc.}
\bvolume{366}
\bpages{1257--1288}.
\bid{doi={10.1090/S0002-9947-2013-05762-5}, issn={0002-9947}, mr={3145731}}
\end{barticle}
%
\bptok{imsref}%
\endbibitem

\bibitem{fop}
\begin{barticle}[mr]
\bauthor{\bsnm{Figari},~\bfnm{R.}\binits{R.}},
\bauthor{\bsnm{Orlandi},~\bfnm{E.}\binits{E.}} \AND
\bauthor{\bsnm{Papanicolaou},~\bfnm{G.}\binits{G.}}
(\byear{1982}).
\btitle{Mean field and {G}aussian approximation for partial differential equations with random coefficients}.
\bjournal{SIAM J. Appl. Math.}
\bvolume{42}
\bpages{1069--1077}.
\bid{doi={10.1137/0142074}, issn={0036-1399}, mr={0673526}}
\end{barticle}
%
\bptok{imsref}%
\endbibitem

\bibitem{fun}
\begin{bincollection}[mr]
\bauthor{\bsnm{Funaki},~\bfnm{Tadahisa}\binits{T.}}
(\byear{2005}).
\btitle{Stochastic interface models}.
In \bbooktitle{Lectures on Probability Theory and Statistics}.
\bseries{Lecture Notes in Math.}
\bvolume{1869}
\bpages{103--274}.
\bpublisher{Springer},
\blocation{Berlin}.
\bid{doi={10.1007/11429579_2}, mr={2228384}}
\end{bincollection}
%
\bptok{imsref}%
\endbibitem

\bibitem{gos}
\begin{barticle}[mr]
\bauthor{\bsnm{Giacomin},~\bfnm{Giambattista}\binits{G.}},
\bauthor{\bsnm{Olla},~\bfnm{Stefano}\binits{S.}} \AND
\bauthor{\bsnm{Spohn},~\bfnm{Herbert}\binits{H.}}
(\byear{2001}).
\btitle{Equilibrium fluctuations for $\nabla\phi$ interface model}.
\bjournal{Ann. Probab.}
\bvolume{29}
\bpages{1138--1172}.
\bid{doi={10.1214/aop/1015345600}, issn={0091-1798}, mr={1872740}}
\end{barticle}
%
\bptok{imsref}%
\endbibitem

\bibitem{glmo}
\begin{barticle}[mr]
\bauthor{\bsnm{Gloria},~\bfnm{Antoine}\binits{A.}} \AND
\bauthor{\bsnm{Mourrat},~\bfnm{Jean-Christophe}\binits{J.-C.}}
(\byear{2012}).
\btitle{Spectral measure and approximation of homogenized coefficients}.
\bjournal{Probab. Theory Related Fields}
\bvolume{154}
\bpages{287--326}.
\bid{doi={10.1007/s00440-011-0370-7}, issn={0178-8051}, mr={2981425}}
\end{barticle}
%
\bptok{imsref}%
\endbibitem

\bibitem{gno-2scale}
\begin{barticle}[mr]
\bauthor{\bsnm{Gloria},~\bfnm{Antoine}\binits{A.}},
\bauthor{\bsnm{Neukamm},~\bfnm{Stefan}\binits{S.}} \AND
\bauthor{\bsnm{Otto},~\bfnm{Felix}\binits{F.}}
(\byear{2014}).
\btitle{An optimal quantitative two-scale expansion in stochastic homogenization of discrete elliptic equations}.
\bjournal{ESAIM Math. Model. Numer. Anal.}
\bvolume{48}
\bpages{325--346}.
\bid{doi={10.1051/m2an/2013110}, issn={0764-583X}, mr={3177848}}
\end{barticle}
%
\bptok{imsref}%
\endbibitem

\bibitem{gno}
\begin{barticle}[mr]
\bauthor{\bsnm{Gloria},~\bfnm{A.}\binits{A.}},
\bauthor{\bsnm{Neukamm},~\bfnm{S.}\binits{S.}} \AND
\bauthor{\bsnm{Otto},~\bfnm{F.}\binits{F.}}
(\byear{2015}).
\btitle{Quantification of ergodicity in stochastic homogenization: Optimal bounds via spectral gap on {G}lauber dynamics}.
\bjournal{Invent. Math.}
\bvolume{199}
\bpages{455--515}.
\bid{doi={10.1007/s00222-014-0518-z}, issn={0020-9910}, mr={3302119}}
\bptnote{check volume, check pages}%
\end{barticle}
%
\bptok{imsref}%
\endbibitem

\bibitem{glotto}
\begin{barticle}[mr]
\bauthor{\bsnm{Gloria},~\bfnm{Antoine}\binits{A.}} \AND
\bauthor{\bsnm{Otto},~\bfnm{Felix}\binits{F.}}
(\byear{2011}).
\btitle{An optimal variance estimate in stochastic homogenization of discrete elliptic equations}.
\bjournal{Ann. Probab.}
\bvolume{39}
\bpages{779--856}.
\bid{doi={10.1214/10-AOP571}, issn={0091-1798}, mr={2789576}}
\end{barticle}
%
\bptok{imsref}%
\endbibitem

\bibitem{glotto2}
\begin{barticle}[mr]
\bauthor{\bsnm{Gloria},~\bfnm{Antoine}\binits{A.}} \AND
\bauthor{\bsnm{Otto},~\bfnm{Felix}\binits{F.}}
(\byear{2012}).
\btitle{An optimal error estimate in stochastic homogenization of discrete elliptic equations}.
\bjournal{Ann. Appl. Probab.}
\bvolume{22}
\bpages{1--28}.
\bid{doi={10.1214/10-AAP745}, issn={1050-5164}, mr={2932541}}
\end{barticle}
%
\bptok{imsref}%
\endbibitem

\bibitem{gubal-1d}
\begin{barticle}[mr]
\bauthor{\bsnm{Gu},~\bfnm{Yu}\binits{Y.}} \AND
\bauthor{\bsnm{Bal},~\bfnm{Guillaume}\binits{G.}}
(\byear{2012}).
\btitle{Random homogenization and convergence to integrals with respect to the {R}osenblatt process}.
\bjournal{J. Differential Equations}
\bvolume{253}
\bpages{1069--1087}.
\bid{doi={10.1016/j.jde.2012.05.007}, issn={0022-0396}, mr={2925905}}
\end{barticle}
%
\bptok{imsref}%
\endbibitem

\bibitem{gubal}
\begin{barticle}[mr]
\bauthor{\bsnm{Gu},~\bfnm{Yu}\binits{Y.}} \AND
\bauthor{\bsnm{Bal},~\bfnm{Guillaume}\binits{G.}}
(\byear{2015}).
\btitle{Fluctuations of parabolic equations with large random potentials}.
\bjournal{Stoch. Partial Differ. Equ. Anal. Comput.}
\bvolume{3}
\bpages{1--51}.
\bid{doi={10.1007/s40072-014-0040-8}, issn={2194-0401}, mr={3312591}}
\bptnote{check volume, check pages}%
\end{barticle}
%
\bptok{imsref}%
\endbibitem

\bibitem{helsjo}
\begin{barticle}[mr]
\bauthor{\bsnm{Helffer},~\bfnm{Bernard}\binits{B.}} \AND
\bauthor{\bsnm{Sj{\"o}strand},~\bfnm{Johannes}\binits{J.}}
(\byear{1994}).
\btitle{On the correlation for {K}ac-like models in the convex case}.
\bjournal{J. Stat. Phys.}
\bvolume{74}
\bpages{349--409}.
\bid{doi={10.1007/BF02186817}, issn={0022-4715}, mr={1257821}}
\end{barticle}
%
\bptok{imsref}%
\endbibitem

\bibitem{kuen}
\begin{barticle}[mr]
\bauthor{\bsnm{K{\"u}nnemann},~\bfnm{Rolf}\binits{R.}}
(\byear{1983}).
\btitle{The diffusion limit for reversible jump processes on {${\bf Z}\sp{d}$} with ergodic random bond conductivities}.
\bjournal{Comm. Math. Phys.}
\bvolume{90}
\bpages{27--68}.
\bid{issn={0010-3616}, mr={0714611}}
\end{barticle}
%
\bptok{imsref}%
\endbibitem

\bibitem{marott}
\begin{barticle}[mr]
\bauthor{\bsnm{Marahrens},~\bfnm{D.}\binits{D.}} \AND
\bauthor{\bsnm{Otto},~\bfnm{F.}\binits{F.}}
(\byear{2016}).
\btitle{Annealed estimates on the Green function}.
\bjournal{Probab. Theory Related Fields}
\bvolume{163}
\bpages{527--573}.
\bid{mr={3418749}}
\end{barticle}
%
\bptok{imsref}%
\endbibitem

\bibitem{mill}
\begin{barticle}[mr]
\bauthor{\bsnm{Miller},~\bfnm{Jason}\binits{J.}}
(\byear{2011}).
\btitle{Fluctuations for the {G}inzburg--{L}andau {$\nabla\phi$} interface model on a bounded domain}.
\bjournal{Comm. Math. Phys.}
\bvolume{308}
\bpages{591--639}.
\bid{doi={10.1007/s00220-011-1315-9}, issn={0010-3616}, mr={2855536}}
\end{barticle}
%
\bptok{imsref}%
\endbibitem

\bibitem{vardecay}
\begin{barticle}[mr]
\bauthor{\bsnm{Mourrat},~\bfnm{Jean-Christophe}\binits{J.-C.}}
(\byear{2011}).
\btitle{Variance decay for functionals of the environment viewed by the particle}.
\bjournal{Ann. Inst. Henri Poincar\'{e} Probab. Stat.}
\bvolume{47}
\bpages{294--327}.
\bid{doi={10.1214/10-AIHP375}, issn={0246-0203}, mr={2779406}}
\end{barticle}
%
\bptok{imsref}%
\endbibitem

\bibitem{homog}
\begin{barticle}[mr]
\bauthor{\bsnm{Mourrat},~\bfnm{Jean-Christophe}\binits{J.-C.}}
(\byear{2014}).
\btitle{Kantorovich distance in the martingale CLT and quantitative homogenization of parabolic equations with random coefficients}.
\bjournal{Probab. Theory Related Fields}
\bvolume{160}
\bpages{279--314}.
\bid{doi={10.1007/s00440-013-0529-5}, issn={0178-8051}, mr={3256815}}
\end{barticle}
%
\bptok{imsref}%
\endbibitem

\bibitem{dl-diff}
\begin{barticle}[mr]
\bauthor{\bsnm{Mourrat},~\bfnm{Jean-Christophe}\binits{J.-C.}}
(\byear{2015}).
\btitle{First-order expansion of homogenized coefficients under {B}ernoulli perturbations}.
\bjournal{J. Math. Pures Appl. (9)}
\bvolume{103}
\bpages{68--101}.
\bid{doi={10.1016/j.matpur.2014.03.008}, issn={0021-7824}, mr={3281948}}
\end{barticle}
%
\bptok{imsref}%
\endbibitem

\bibitem{fluct2}
\begin{bmisc}[auto:parserefs-M02]
\bauthor{\bsnm{Mourrat},~\bfnm{J.-C.}\binits{J.-C.}} \AND
\bauthor{\bsnm{Gu},~\bfnm{Y.}\binits{Y.}}
(\byear{2015}).
\bhowpublished{Scaling limit of fluctuations in stochastic homogenization.
Unpublished manuscript. Available at \arxivurl{arXiv:1503.00578}.}
\end{bmisc}
%
\bptok{imsref}%
\endbibitem

\bibitem{fluct}
\begin{bmisc}[auto:parserefs-M02]
\bauthor{\bsnm{Mourrat},~\bfnm{J.-C.}\binits{J.-C.}} \AND
\bauthor{\bsnm{Nolen},~\bfnm{J.}\binits{J.}}
(\byear{2015}).
\bhowpublished{Scaling limit of the corrector in stochastic homogenization.
Unpublished manuscript. Available at \arxivurl{arXiv:1502.07440}.}
\end{bmisc}
%
\bptok{imsref}%
\endbibitem

\bibitem{nadspe}
\begin{barticle}[mr]
\bauthor{\bsnm{Naddaf},~\bfnm{Ali}\binits{A.}} \AND
\bauthor{\bsnm{Spencer},~\bfnm{Thomas}\binits{T.}}
(\byear{1997}).
\btitle{On homogenization and scaling limit of some gradient perturbations of a massless free field}.
\bjournal{Comm. Math. Phys.}
\bvolume{183}
\bpages{55--84}.
\bid{doi={10.1007/BF02509796}, issn={0010-3616}, mr={1461951}}
\end{barticle}
%
\bptok{imsref}%
\endbibitem

\bibitem{nadspe-unpub}
\begin{bmisc}[auto:parserefs-M02]
\bauthor{\bsnm{Naddaf},~\bfnm{A.}\binits{A.}} \AND
\bauthor{\bsnm{Spencer},~\bfnm{T.}\binits{T.}}
(\byear{1998}).
\bhowpublished{Estimates on the variance of some homogenization problems.
Unpublished manuscript.}
\end{bmisc}
%
\bptok{imsref}%
\endbibitem

\bibitem{papvar}
\begin{bincollection}[mr]
\bauthor{\bsnm{Papanicolaou},~\bfnm{G.~C.}\binits{G.~C.}} \AND
\bauthor{\bsnm{Varadhan},~\bfnm{S.~R.~S.}\binits{S.~R.~S.}}
(\byear{1981}).
\btitle{Boundary value problems with rapidly oscillating random coefficients}.
In \bbooktitle{Random Fields, {V}ol. I, II ({E}sztergom, 1979)}.
\bseries{Colloquia Mathematica Societatis J\'anos Bolyai}
\bvolume{27}
\bpages{835--873}.
\bpublisher{North-Holland},
\blocation{Amsterdam}.
\bid{mr={0712714}}
\end{bincollection}
%
\bptok{imsref}%
\endbibitem

\bibitem{sjo}
\begin{barticle}[mr]
\bauthor{\bsnm{Sj{\"o}strand},~\bfnm{J.}\binits{J.}}
(\byear{1996}).
\btitle{Correlation asymptotics and {W}itten {L}aplacians}.
\bjournal{Algebra i Analiz}
\bvolume{8}
\bpages{160--191}.
\bid{issn={0234-0852}, mr={1392018}}
\end{barticle}
%
\bptok{imsref}%
\endbibitem

\bibitem{yuri}
\begin{barticle}[mr]
\bauthor{\bsnm{Yurinski{\u\i}},~\bfnm{V.~V.}\binits{V.~V.}}
(\byear{1986}).
\btitle{Averaging of symmetric diffusion in a random medium}.
\bjournal{Sibirsk. Mat. Zh.}
\bvolume{27}
\bpages{167--180, 215}.
\bid{issn={0037-4474}, mr={0867870}}
\bptnote{check pages}%
\end{barticle}
%
\bptok{imsref}%
\endbibitem

\end{thebibliography}
\end{document}